\documentclass[12pt]{article}
\usepackage{wasysym,amssymb,indentfirst,cite,bibspacing,mathrsfs}
\begin{document}
\baselineskip 16pt
\title{\textbf{On the $\pi$$\mathfrak{F}$-norm and the $\mathfrak{H}$-$\mathfrak{F}$-norm of a finite group}\thanks{Research is supported by a NNSF grant of China (grant \#11371335) and Research Fund for the Doctoral Program of Higher Education of China (Grant 20113402110036).}}
\author{Xiaoyu Chen, Wenbin Guo\thanks{Corresponding author.}\\
{\small  Department of Mathematics, University of Science and Technology of China,}\\ {\small Hefei 230026, P. R. China}\\
 {\small E-mails: jelly@mail.ustc.edu.cn, $\,$wbguo@ustc.edu.cn}
}
\date{}
\maketitle
\begin{abstract}
Let $\mathfrak{H}$ be a Fitting class and $\mathfrak{F}$ a formation. We call a subgroup $\mathcal{N}_{\mathfrak{H},\mathfrak{F}}(G)$ of a finite group $G$ the $\mathfrak{H}$-$\mathfrak{F}$-norm of $G$ if $\mathcal{N}_{\mathfrak{H},\mathfrak{F}}(G)$ is the intersection of the normalizers of the products of the $\mathfrak{F}$-residuals of all subgroups of $G$ and the $\mathfrak{H}$-radical of $G$. Let $\pi$ denote a set of primes and let $\mathfrak{G}_\pi$ denote the class of all finite $\pi$-groups. We call the subgroup $\mathcal{N}_{\mathfrak{G}_\pi,\mathfrak{F}}(G)$ of $G$ the $\pi\mathfrak{F}$-norm of $G$. A normal subgroup $N$ of $G$ is called $\pi\mathfrak{F}$-hypercentral in $G$ if either $N=1$ or $N>1$ and every $G$-chief factor below $N$ of order divisible by at least one prime in $\pi$ is $\mathfrak{F}$-central in $G$. Let $Z_{\pi\mathfrak{F}}(G)$ denote the $\pi\mathfrak{F}$-hypercentre of $G$, that is, the product of all $\pi\mathfrak{F}$-hypercentral normal subgroups of $G$. In this paper, we study the properties of the $\mathfrak{H}$-$\mathfrak{F}$-norm, especially of the $\pi\mathfrak{F}$-norm of a finite group $G$. In particular, we investigate the relationship between the $\pi'\mathfrak{F}$-norm and the $\pi\mathfrak{F}$-hypercentre of $G$.\par
\end{abstract}
\renewcommand{\thefootnote}{\empty}
\footnotetext{Keywords: norm, $\pi\mathfrak{F}$-norm, $\pi\mathfrak{F}$-hypercentre, $\mathfrak{F}$-maximal subgroups, $\mathfrak{F}$-critical group.}
\footnotetext{Mathematics Subject Classification (2000): 20D10, 20D15.}

\section{Introduction}
\noindent All groups considered in this paper are finite, and all classes of groups $\mathfrak{X}$ mentioned are non-empty. $G$ always denotes a group, $p$ denotes a prime, $\pi$ denotes a set of primes, and $\mathbb{P}$ denotes the set of all primes. Also, let $\pi(G)$ denote the set of all prime divisors of the order of $G$, and let $\pi(\mathfrak{X})=\bigcup\{\pi(G):G\in \mathfrak{X}\}$ for a class of groups $\mathfrak{X}$.\par
Recall that a class of groups $\mathfrak{F}$ is called a formation if $\mathfrak{F}$ is closed under taking homomorphic images and subdirect products. A formation $\mathfrak{F}$ is said to be saturated if $G\in \mathfrak{F}$ whenever $G/\Phi(G)\in \mathfrak{F}$. The $\mathfrak{F}$-residual of $G$, denoted by $G^\mathfrak{F}$, is the smallest normal subgroup $N$ of $G$ with $G/N\in \mathfrak{F}$. The formation product $\mathfrak{X}\circ\mathfrak{F}$ of a class of groups $\mathfrak{X}$ and a formation $\mathfrak{F}$ is the class of all groups $G$ such that $G^\mathfrak{F}\in\mathfrak{X}$.
A class of groups $\mathfrak{H}$ is called a Fitting class if $\mathfrak{H}$ is closed under taking normal subgroups and products of normal $\mathfrak{H}$-subgroups. The $\mathfrak{H}$-radical of $G$, denoted by $G_\mathfrak{H}$, is the maximal normal $\mathfrak{H}$-subgroup of $G$. The Fitting product $\mathfrak{H}\diamond\mathfrak{X}$ of a Fitting class $\mathfrak{H}$ and a class of groups $\mathfrak{X}$ is the class of all groups $G$ such that $G/G_\mathfrak{H}\in \mathfrak{X}$. A class of groups $\mathfrak{B}$ is called a Fitting formation if $\mathfrak{B}$ is both a formation and a Fitting class. Note that for a Fitting formation $\mathfrak{B}$, a formation $\mathfrak{F}$ and a Fitting class $\mathfrak{H}$, $\mathfrak{H}\diamond(\mathfrak{B}\circ\mathfrak{F})=(\mathfrak{H}\diamond\mathfrak{B})\circ\mathfrak{F}$ always holds, and we denote it by $\mathfrak{H}\diamond\mathfrak{B}\circ\mathfrak{F}$.\par
The class of the groups of order $1$ is denoted by $1$, and the class of all finite groups is denoted by $\mathfrak{G}$. We use $\mathfrak{S}$ (resp. $\mathfrak{N}$, $\mathfrak{U}$, $\mathfrak{A}$) to denote the class of finite solvable (resp. nilpotent, supersolvable, abelian) groups and $\mathfrak{S}_\pi$ (resp. $\mathfrak{N}_\pi$, $\mathfrak{U}_\pi$) to denote the class of finite $\pi$-solvable (resp. $\pi$-nilpotent, $\pi$-supersolvable) groups. Also, the symbol $\mathfrak{G}_\pi$ denotes the class of all finite $\pi$-groups.\par
A formation function $f$ is a local function $f$: $\mathbb{P}\rightarrow\{\mbox{classes of groups}\}$ such that $f(p)$ is a formation for all $p\in\mathbb{P}$. Let $LF(f)$ denote the set of all groups $G$ whose chief factors $L/K$ are all $f$-central in $G$, that is, $G/C_G(L/K)\in f(p)$ for all $p\in \pi(L/K)$.
The canonical local definition of a saturated formation $\mathfrak{F}$ is the uniquely determined formation function $F$ such that $\mathfrak{F}=LF(F)$, $F(p)\subseteq \mathfrak{F}$ and $\mathfrak{G}_p\circ F(p)=F(p)$ for all $p\in\mathbb{P}$ (for details, see \cite[Chap. IV]{Doe}).\par
Following \cite[Chap. II]{Doe}, for a class of groups $\mathfrak{X}$, we define closure operations as follows: $\mathtt{S}\mathfrak{X}=(G:G\leq H\ \mbox{for some}\ H\in \mathfrak{X})$; $\mathtt{S_n}\mathfrak{X}=(G:G\ \mbox{is subnormal in }H\ \mbox{for some}\ H\in \mathfrak{X})$; $\mathtt{Q}\mathfrak{X}=(G: \mbox{ there exist}\ H\in \mathfrak{X}$ $\mbox{and an epimorphism from}$ $H\ \mbox{onto}\ G)$; $\mathtt{E}\mathfrak{X}=(G:\mbox{there exists a series}$ $\mbox{of subgroups of }G:$ $1=G_0\unlhd G_1\unlhd\cdots\unlhd G_n=G$ $\mbox{with each }G_i/G_{i-1}\in \mathfrak{X})=\bigcup^{\infty}_{r=1}\mathfrak{X}^r$.\par
Recall that the norm $\mathcal{N}(G)$ of $G$ is the intersection of the normalizers of all subgroups of $G$, and the Wielandt subgroup $\omega(G)$ of $G$ is the intersection of the normalizers of all subnormal subgroups of $G$. These concepts were introduced by R. Baer \cite{Bae} and H. Wielandt \cite{Wie} in 1934 and 1958, respectively. Much investigation has focused on using the concepts of the norm and the Wielandt subgroup to determine the structure of finite groups (see, for example,\cite{Sch,Bae1,Bae2,Bei,Bry,Cam,Cos,Kap,Orm}).\par
Recently, Li and Shen \cite{Li} considered the intersection of the normalizers of the derived subgroups of all subgroups of $G$. Also, in \cite{Gon} and \cite{She}, the authors considered the intersection of the normalizers of the nilpotent residuals of all subgroups of $G$. Furthermore, for a formation $\mathfrak{F}$, Su and Wang \cite{Su} investigated the intersection of the normalizers of the $\mathfrak{F}$-residuals of all subgroups of $G$ and the intersection of the normalizers of the products of the $\mathfrak{F}$-residuals of all subgroups of $G$ and $O_{p'}(G)$. As a continuation of the above ideas, we now introduce the notion of $\mathfrak{H}$-$\mathfrak{F}$-norm as follows:\par
\medskip
\noindent\textbf{Definition 1.1.} Let $\mathfrak{H}$ be a Fitting class and $\mathfrak{F}$ a formation. We call a subgroup $\mathcal{N}_{\mathfrak{H},\mathfrak{F}}(G)$ of $G$ the \textit{$\mathfrak{H}$-$\mathfrak{F}$-norm} of $G$ if $\mathcal{N}_{\mathfrak{H},\mathfrak{F}}(G)$ is the intersection of the normalizers of the products of the $\mathfrak{F}$-residuals of all subgroups of $G$ and the $\mathfrak{H}$-radical of $G$, that is, $$\mathcal{N}_{\mathfrak{H},\mathfrak{F}}(G)=\bigcap_{H\leq G}N_G(H^{\mathfrak{F}}G_{\mathfrak{H}}).$$ In particular, when $\mathfrak{H}=1$, the subgroup $\mathcal{N}_{1,\mathfrak{F}}(G)$ of $G$ is called the \textit{$\mathfrak{F}$-norm} of $G$, and we denote it by $\mathcal{N}_{\mathfrak{F}}(G)$, that is, $$\mathcal{N}_{\mathfrak{F}}(G)=\bigcap_{H\leq G}N_G(H^\mathfrak{F});$$ when $\mathfrak{H}=\mathfrak{G}_\pi$, the subgroup $\mathcal{N}_{\mathfrak{G}_\pi,\mathfrak{F}}(G)$ of $G$ is called the \textit{$\pi\mathfrak{F}$-norm} of $G$, and we denote it by $\mathcal{N}_{\pi\mathfrak{F}}(G)$, that is, $$\mathcal{N}_{\pi\mathfrak{F}}(G)=\bigcap_{H\leq G}N_G(H^\mathfrak{F}O_\pi(G)).$$\par
\medskip
\noindent\textbf{Definition 1.2.} Let $\mathcal{N}_{\mathfrak{H},\mathfrak{F}}^{0}(G)=1$ and $\mathcal{N}_{\mathfrak{H},\mathfrak{F}}^{i}(G)/\mathcal{N}_{\mathfrak{H},\mathfrak{F}}^{i-1}(G)=\mathcal{N}_{\mathfrak{H},\mathfrak{F}}(G/\mathcal{N}_{\mathfrak{H},\mathfrak{F}}^{i-1}(G))$ for $i=1,2,\cdots$. Then there exists a series of subgroups of $G$: $$1=\mathcal{N}_{\mathfrak{H},\mathfrak{F}}^{0}(G)\leq \mathcal{N}_{\mathfrak{H},\mathfrak{F}}^{1}(G) \leq \mathcal{N}_{\mathfrak{H},\mathfrak{F}}^{2}(G)\cdots \leq \mathcal{N}_{\mathfrak{H},\mathfrak{F}}^{n}(G)=\mathcal{N}_{\mathfrak{H},\mathfrak{F}}^{n+1}(G)=\cdots.$$ Denote $\mathcal{N}_{\mathfrak{H},\mathfrak{F}}^{\infty}(G)$ the terminal term of this ascending series. In particular, when $\mathfrak{H}=1$, we denote $\mathcal{N}_{1,\mathfrak{F}}^{\infty}(G)$ by $\mathcal{N}_{\mathfrak{F}}^{\infty}(G)$; when $\mathfrak{H}=\mathfrak{G}_\pi$, we denote $\mathcal{N}_{\mathfrak{G}_\pi,\mathfrak{F}}^{\infty}(G)$ by $\mathcal{N}_{\pi\mathfrak{F}}^{\infty}(G)$.\par
\medskip
Let $\mathfrak{F}$ be a formation. A $G$-chief factor $L/K$ is said to be $\mathfrak{F}$-central in $G$ if $(L/K)\rtimes (G/C_G(L/K))\in \mathfrak{F}$. Following \cite{Guo}, a normal subgroup $N$ of $G$ is called \textit{$\pi\mathfrak{F}$-hypercentral} in $G$ if either $N=1$ or $N>1$ and every $G$-chief factor below $N$ of order divisible by at least one prime in $\pi$ is $\mathfrak{F}$-central in $G$. Let $Z_{\pi\mathfrak{F}}(G)$ denote the \textit{$\pi\mathfrak{F}$-hypercentre} of $G$, that is, the product of all $\pi\mathfrak{F}$-hypercentral normal subgroups of $G$. The $\mathbb{P}\mathfrak{F}$-hypercentre of $G$ is called the $\mathfrak{F}$-hypercentre of $G$, and we denote it by $Z_{\mathfrak{F}}(G)$.\par
Let $\mathfrak{X}$ be a class of groups. Recall that a subgroup $U$ of $G$ is called \textit{$\mathfrak{X}$-maximal} in $G$ if $U\in \mathfrak{X}$ and $G$ does not have a subgroup $V$ such that $U<V$ and $V\in \mathfrak{X}$. Following \cite{Ski}, we use $\mbox{Int}_\mathfrak{X}(G)$ to denote the intersection of all $\mathfrak{X}$-maximal subgroups of $G$.\par
In \cite[Remark 4]{Bei1}, J. C. Beidleman and H. Heineken observed that $\mathcal{N}_{\mathfrak{N}_c}^{\infty}(G)$ coincides with $\mbox{Int}_{\mathfrak{N}\circ\mathfrak{N}_c}(G)$ for every group $G$, where $\mathfrak{N}_c$ denotes the class of nilpotent groups of class at most $c$. In \cite{Ski}, A. N. Skiba gave conditions under which the $\mathfrak{F}$-hypercentre $Z_\mathfrak{F}(G)$ coincides with $\mbox{Int}_{\mathfrak{F}}(G)$ for every group $G$. Also, Guo and A. N. Skiba \cite{Guo} gave conditions under which the $\pi\mathfrak{F}$-hypercentre $Z_{\pi\mathfrak{F}}(G)$ coincides with $\mbox{Int}_{\mathfrak{F}}(G)$ for every group $G$.\par
Motivated by the above observations, the following questions naturally arise:\par
\medskip
\noindent\textbf{Problem (I).} \textit{Under what conditions $\mathcal{N}_{\mathfrak{F}}^{\infty}(G)$ coincides with the $\mathfrak{N}\circ\mathfrak{F}$-hypercentre $Z_{\mathfrak{N}\circ\mathfrak{F}}(G)$\textup{?} More generally, under what conditions $\mathcal{N}_{\pi'\mathfrak{F}}^{\infty}(G)$ coincides with the $\pi(\mathfrak{N}\circ\mathfrak{F})$-hypercentre $Z_{\pi({\mathfrak{N}\circ\mathfrak{F}})}(G)$\textup{?}}\par
\medskip
\noindent\textbf{Problem (II).} \textit{Under what conditions $\mathcal{N}_{\mathfrak{F}}^{\infty}(G)$ coincides with $\mbox{\textup{Int}}_{\mathfrak{N}\circ\mathfrak{F}}(G)$\textup{?} More generally, under what conditions $\mathcal{N}_{\pi'\mathfrak{F}}^{\infty}(G)$ coincides with $\mbox{\textup{Int}}_{\mathfrak{N}_\pi\circ\mathfrak{F}}(G)$}?\par
\medskip
For a class of groups $\mathfrak{X}$, a group $G$ is called $\mathtt{S}$-critical for $\mathfrak{X}$ if $G\notin \mathfrak{X}$ but all proper subgroups of $G$ belong to $\mathfrak{X}$. Let $\mbox{Crit}_\mathtt{S}(\mathfrak{X})$ denote the set of all groups $G$ which are $\mathtt{S}$-critical for $\mathfrak{X}$. For convenience of statement, we give the following definition.\par
\medskip
\noindent\textbf{Definition 1.3.} We say that a formation $\mathfrak{F}$ satisfies:\par
(1) \textit{The $\pi$-boundary condition} (I) if $\mbox{Crit}_\mathtt{S}(\mathfrak{F})\subseteq \mathfrak{N_{\pi}}\circ\mathfrak{F}$ (equivalently, $\mbox{Crit}_\mathtt{S}(\mathfrak{F})\subseteq \mathfrak{S_{\pi}}\circ\mathfrak{F}$, see Lemma 2.7 below).\par
(2) \textit{The $\pi$-boundary condition} (II) if for any $p\in \pi$, $\mbox{Crit}_\mathtt{S}(\mathfrak{G}_p\circ\mathfrak{F})\subseteq \mathfrak{S}_\pi\circ\mathfrak{F}$.\par
(3) \textit{The $\pi$-boundary condition} (III) if for any $p\in \pi$, $\mbox{Crit}_\mathtt{S}(\mathfrak{G}_p\circ\mathfrak{F})\subseteq \mathfrak{N}_\pi\circ\mathfrak{F}$.\par
(4) \textit{The $\pi$-boundary condition \textup{(III)} in $\mathfrak{S}$} if for any $p\in \pi$, $\mbox{Crit}_\mathtt{S}(\mathfrak{G}_p\circ\mathfrak{F})\cap \mathfrak{S}\subseteq \mathfrak{N}_\pi\circ\mathfrak{F}$.\par
\medskip
Note that a formation $\mathfrak{F}$ satisfies the $\pi$-boundary condition (III) (resp. the $\pi$-boundary condition \textup{(III)} in $\mathfrak{S}$) if and only if $\mathfrak{N}_\pi\circ\mathfrak{F}$ satisfies the $\pi$-boundary condition (resp. the $\pi$-boundary condition in $\mathfrak{S}$) in the sense of \cite{Guo}.\par
\medskip
\noindent\textbf{Remark 1.4.} If a formation $\mathfrak{F}$ satisfies the $\pi$-boundary condition (II), then clearly, $\mathfrak{F}$ satisfies the $\pi$-boundary condition (I). However, the converse does not hold. For example, let $\pi=\mathbb{P}$ and $\mathfrak{F}=\mathfrak{N}_3$. By \cite[Chap. IV, Satz 5.4]{Hup}, $\mbox{Crit}_\mathtt{S}(\mathfrak{N}_3)\subseteq \mathfrak{N}\circ\mathfrak{N}_3$. Now let $G=A_5$, where $A_5$ is the alternating group of degree 5. Then $G\in \mbox{Crit}_\mathtt{S}(\mathfrak{G}_3\circ\mathfrak{N}_3)$, but $G\notin \mathfrak{S}\circ\mathfrak{N}_3$. Hence $\mbox{Crit}_\mathtt{S}(\mathfrak{G}_3\circ\mathfrak{N}_3)\nsubseteq \mathfrak{S}\circ\mathfrak{N}_3$.\par
\medskip
\noindent\textbf{Remark 1.5.} If a formation $\mathfrak{F}$ satisfies the $\pi$-boundary condition (III), then $\mathfrak{F}$ satisfies the $\pi$-boundary condition (II). However, the converse does not hold. For example, let $\pi=\mathbb{P}$ and $\mathfrak{F}=\mathfrak{G}_3$. For any prime $p\neq 3$, $\mbox{Crit}_\mathtt{S}(\mathfrak{G}_p\circ\mathfrak{G}_3)\subseteq \mathfrak{N}_3\cup \mbox{Crit}_\mathtt{S}(\mathfrak{N}_3)$. If there exists a group $H$ such that $H\in \mbox{Crit}_\mathtt{S}(\mathfrak{G}_p\circ\mathfrak{G}_3)\backslash (\mathfrak{S}\circ\mathfrak{G}_3)$, then by \cite[Chap. IV, Satz 5.4]{Hup}, we have that $H\in \mathfrak{N}_3$. Hence $H$ has the normal 3-complement $A$. If $A<H$, then $A\in \mathfrak{G}_p\circ\mathfrak{G}_3\subseteq \mathfrak{S}$, and thereby $H\in \mathfrak{S}$, a contradiction. Therefore, $H=A\in \mathfrak{G}_p\cup\mbox{Crit}_\mathtt{S}(\mathfrak{G}_p)\subseteq \mathfrak{S}$, also a contradiction. This shows that $\mbox{Crit}_\mathtt{S}(\mathfrak{G}_p\circ\mathfrak{G}_3)\subseteq \mathfrak{S}\circ\mathfrak{G}_3$, and so $\mathfrak{G}_3$ satisfies the $\mathbb{P}$-boundary condition (II). Now let $G=S_3$, where $S_3$ is the symmetric group of degree 3. Then it is easy to see that $G\in \mbox{Crit}_\mathtt{S}(\mathfrak{G}_2\circ\mathfrak{G}_3)$, but $G\notin \mathfrak{N}\circ\mathfrak{G}_3$. Hence $\mbox{Crit}_\mathtt{S}(\mathfrak{G}_2\circ\mathfrak{G}_3)\nsubseteq \mathfrak{N}\circ\mathfrak{G}_3$.\par
\medskip
Firstly, we give a characterization of $\mathfrak{H}\diamond\mathfrak{N}\circ\mathfrak{F}$-groups by using their $\mathfrak{H}$-$\mathfrak{F}$-norms.\par
\medskip
\noindent\textbf{Theorem A.} Let $\mathfrak{H}$ be a saturated Fitting formation such that $\mathfrak{G}_{\pi'}\subseteq\mathfrak{H}=\mathtt{E}\mathfrak{H}$ and $\mathfrak{F}$ a formation such that $\mathfrak{F}=\mathtt{S}\mathfrak{F}$. Suppose that one of the following holds:\par
(i) $G^{\mathfrak{H}\diamond\mathfrak{N}\circ\mathfrak{F}}\in \mathfrak{S_{\pi}}$.\par
(ii) $\mathfrak{F}$ satisfies the $\pi$-boundary condition (I).\par
\noindent Then the following statements are equivalent:\par
(1) $G\in \mathfrak{H}\diamond\mathfrak{N}\circ\mathfrak{F}$.\par
(2) $G/\mathcal{N}_{\mathfrak{H},\mathfrak{F}}(G)\in \mathfrak{H}\diamond\mathfrak{N}\circ\mathfrak{F}$.\par
(3) $G/\mathcal{N}_{\mathfrak{H},\mathfrak{F}}^{\infty}(G)\in \mathfrak{H}\diamond\mathfrak{N}\circ\mathfrak{F}$.\par
(4) $\mathcal{N}_{\mathfrak{H},\mathfrak{F}}(G/N)>1$ for every proper normal subgroup $N$ of $G$.\par
(5) $G=\mathcal{N}_{\mathfrak{H},\mathfrak{F}}^{\infty}(G)$.\par
\medskip
The main purpose of this paper is to give answers to Problem (I) and (II). In the universe of all groups, we prove:\par
\medskip
\noindent\textbf{Theorem B.} Let $\mathfrak{F}$ be a formation such that $\mathfrak{F}=\mathtt{S}\mathfrak{F}$. Then:\par
(1) If $\mathfrak{F}$ satisfies the $\pi$-boundary condition (II), then $\mathcal{N}_{\pi'\mathfrak{F}}^{\infty}(G)=Z_{\pi(\mathfrak{N}\circ\mathfrak{F})}(G)$ holds for every group $G$.\par
(2) If $\mathcal{N}_{\pi'\mathfrak{F}}^{\infty}(G)=Z_{\pi(\mathfrak{N}\circ\mathfrak{F})}(G)$ holds for every group $G$, then $\mathfrak{F}$ satisfies the $\pi$-boundary condition (I).\par
(3) $\mathcal{N}_{\pi'\mathfrak{F}}^{\infty}(G)=Z_{\pi(\mathfrak{N}\circ\mathfrak{F})}(G)$ holds for every group $G$ if and only if $\mathcal{N}_{\pi'\mathfrak{F}}^{\infty}(G)=Z_{\pi(\mathfrak{N}\circ\mathfrak{F})}(G)$ holds for every group $G\in \bigcup_{p\in\pi} (\mbox{Crit}_\mathtt{S}(\mathfrak{G}_p\circ\mathfrak{F})\backslash (\mathfrak{S}_\pi\circ\mathfrak{F}))$.\par
\medskip
\noindent\textbf{Remark 1.6.} The converse of statement (2) of Theorem B does not hold. For example, let $\pi=\mathbb{P}$ and $\mathfrak{F}=\mathfrak{U}$. By K. Doerk's result \cite{Doe1}, $\mbox{Crit}_\mathtt{S}(\mathfrak{U})\subseteq \mathfrak{N}\circ\mathfrak{U}$. This means that $\mathfrak{U}$ satisfies the $\mathbb{P}$-boundary condition (I). Let $A$ be the $2$-Frattini module of $A_5$, where $A_5$ is the alternating group of degree 5. By \cite[Example 1]{Gri}, the dimension of $A$ is 5. Then by \cite[Appendix $\beta$, Proposition $\beta.5$]{Doe}, there exists a Frattini extension $G$ such that $G/A\cong A_5$ and $A=\Phi(G)$. Now we show that $\mathcal{N}_{\mathfrak{U}}(G)=\Phi(G)$. As $\mathcal{N}_{\mathfrak{U}}(G)<G$, it will suffice to prove that for any subgroup $H$ of $G$, $\Phi(G)\leq N_G(H^\mathfrak{U})$. If $H/H\cap \Phi(G)\in \mathfrak{U}$, then $H^\mathfrak{U}\leq \Phi(G)$, and so $\Phi(G)\leq N_G(H^\mathfrak{U})$. Hence, consider that $H/H\cap \Phi(G)\notin \mathfrak{U}$. Since $G/\Phi(G)\cong A_5$, $H\Phi(G)/\Phi(G)\cong A_4$, where $A_4$ is the alternating group of degree 4. This implies that $H\Phi(G)$ is a Hall $5'$-subgroup of $G$, and thereby $H$ is a Hall $5'$-subgroup of $G$. Thus $\Phi(G)\leq H$, and consequently $\Phi(G)\leq N_G(H^\mathfrak{U})$. Therefore, $\mathcal{N}_{\mathfrak{U}}(G)=\Phi(G)$. If $\mathcal{N}_{\mathfrak{U}}^{\infty}(G)=Z_{\mathfrak{N}\circ\mathfrak{U}}(G)$, then $Z_{\mathfrak{N}\circ\mathfrak{U}}(G)=\Phi(G)$. Since $G^{\mathfrak{N}\circ\mathfrak{U}}=G$, by \cite[Chap. IV, Theorem 6.10]{Doe}, $Z(G)=\Phi(G)$. It follows that $G$ is quasisimple. By \cite[Table 4.1]{Gor}, the Schur multiplier of $A_5$ is a cyclic group of order 2, a contradiction. Hence $\mathcal{N}_{\mathfrak{U}}^{\infty}(G)\neq Z_{\mathfrak{N}\circ\mathfrak{U}}(G)$. Besides, we currently do not know whether the converse of statement (1) of Theorem B is true or not.\par
\medskip
\noindent\textbf{Theorem C.} Let $\mathfrak{F}$ be a formation such that $\mathfrak{F}=\mathtt{S}\mathfrak{F}$. Then the following statements are equivalent:\par
(1) $\mathcal{N}_{\pi'\mathfrak{F}}^{\infty}(G)=\mbox{Int}_{\mathfrak{N}_\pi\circ\mathfrak{F}}(G)$ holds for every group $G$.\par
(2) $Z_{\pi(\mathfrak{N}\circ\mathfrak{F})}(G)=\mbox{Int}_{\mathfrak{N}_\pi\circ\mathfrak{F}}(G)$ holds for every group $G$.\par
(3) $\mathfrak{F}$ satisfies the $\pi$-boundary condition (III).\par
\medskip
In the universe of all solvable groups, we prove:\par
\medskip
\noindent\textbf{Theorem D.} Let $\mathfrak{F}$ be a formation such that $\mathfrak{F}=\mathtt{S}\mathfrak{F}$. Then $\mathcal{N}_{\pi'\mathfrak{F}}^{\infty}(G)=Z_{\pi(\mathfrak{N}\circ\mathfrak{F})}(G)$ holds for every group $G\in \mathfrak{S_{\pi}}\circ\mathfrak{F}$.\par
\medskip
\noindent\textbf{Theorem E.} Let $\mathfrak{F}$ be a formation such that $\mathfrak{F}=\mathtt{S}\mathfrak{F}$. Then the following statements are equivalent:\par
(1) $\mathcal{N}_{\pi'\mathfrak{F}}^{\infty}(G)=\mbox{Int}_{\mathfrak{N}_\pi\circ\mathfrak{F}}(G)$ holds for every $G\in \mathfrak{S}$.\par
(2) $Z_{\pi(\mathfrak{N}\circ\mathfrak{F})}(G)=\mbox{Int}_{\mathfrak{N}_\pi\circ\mathfrak{F}}(G)$ holds for every $G\in \mathfrak{S}$.\par
(3) $\mathfrak{F}$ satisfies the $\pi$-boundary condition (III) in $\mathfrak{S}$.\par
\section{Preliminaries}
The following two lemmas are well known.\par
\medskip
\noindent\textbf{Lemma 2.1.} Let $\mathfrak{F}$ be a formation. Suppose that $H\leq G$ and $N\unlhd G$. Then:\par
(1) $G^\mathfrak{F}N/N=(G/N)^\mathfrak{F}$.\par
(2) If $\mathfrak{F}=\mathtt{S}\mathfrak{F}$ (resp. $\mathfrak{F}=\mathtt{S_n}\mathfrak{F}$), then $H^\mathfrak{F}\leq G^\mathfrak{F}\cap H$ (resp. $N^\mathfrak{F}\leq G^\mathfrak{F}\cap N$).\par
\medskip
\noindent\textbf{Lemma 2.2.} Let $\mathfrak{H}$ be a Fitting class. Suppose that $H\leq G$ and $N\unlhd G$. Then:\par
(1) $G_\mathfrak{H}\cap N=N_\mathfrak{H}$.\par
(2) If $\mathfrak{H}=\mathtt{S}\mathfrak{H}$, then $G_\mathfrak{H}\cap H\leq H_\mathfrak{H}$.\par
(3) If $\mathfrak{H}=\mathtt{Q}\mathfrak{H}$, then $G_\mathfrak{H}N/N\leq (G/N)_\mathfrak{H}$.\par
(4) If $\mathfrak{H}=\mathtt{E}\mathfrak{H}$ and $N\leq G_\mathfrak{H}$, then $(G/N)_\mathfrak{H}\leq G_\mathfrak{H}/N$.\par
\medskip
\noindent\textbf{Lemma 2.3.} Let $\mathfrak{H}$ be a Fitting class and $\mathfrak{F}$ a formation. Suppose that $H\leq G$ and $N\unlhd G$. Then:\par
(1) $\mathcal{N}_{\mathfrak{H},\mathfrak{F}}(G)\cap N\leq \mathcal{N}_{\mathfrak{H},\mathfrak{F}}(N)$.\par
(2) If $\mathfrak{H}=\mathtt{S}\mathfrak{H}$, then $\mathcal{N}_{\mathfrak{H},\mathfrak{F}}(G)\cap H\leq \mathcal{N}_{\mathfrak{H},\mathfrak{F}}(H)$.\par
(3) If $\mathfrak{H}=\mathtt{Q}\mathfrak{H}$, then $\mathcal{N}_{\mathfrak{H},\mathfrak{F}}(G)N/N\leq \mathcal{N}_{\mathfrak{H},\mathfrak{F}}(G/N)$.\par
(4) If $\mathfrak{F}=\mathtt{S}\mathfrak{F}$ and $G\in \mathfrak{H}\diamond\mathfrak{N}\circ\mathfrak{F}$, then either $G=1$ or $\mathcal{N}_{\mathfrak{H},\mathfrak{F}}(G)>1$.\par
\medskip
\noindent\textit{Proof.} (1) By definition and Lemma 2.2(1), $\mathcal{N}_{\mathfrak{H},\mathfrak{F}}(G)\cap N=(\bigcap_{H\leq G}N_G(H^{\mathfrak{F}}G_{\mathfrak{H}}))\cap N\leq\bigcap_{H\leq N}\linebreak[4]N_N(H^{\mathfrak{F}}G_{\mathfrak{H}})\leq \bigcap_{H\leq N}N_N(H^{\mathfrak{F}}(G_{\mathfrak{H}}\cap N))=\bigcap_{H\leq N}N_N(H^{\mathfrak{F}}N_{\mathfrak{H}})=\mathcal{N}_{\mathfrak{H},\mathfrak{F}}(N)$.\par
The proof of statement (2) is similar to (1).\par
(3) By definition, Lemma 2.1(1) and Lemma 2.2(3), $\mathcal{N}_{\mathfrak{H},\mathfrak{F}}(G)N/N=$$(\bigcap_{H\leq G}N_G(H^{\mathfrak{F}}G_{\mathfrak{H}}))N/N\linebreak[4]\leq \bigcap_{N\leq H\leq G}$$N_G(H^{\mathfrak{F}}$$G_{\mathfrak{H}})/N$
$\leq \bigcap_{H/N\leq G/N}N_{G/N}((H^{\mathfrak{F}}N/N)$$(G/N)_{\mathfrak{H}})=\mathcal{N}_{\mathfrak{H},\mathfrak{F}}(G/N)$.\par
(4) We may suppose that $G>1$ and $G_\mathfrak{H}=1$. Since $G\in \mathfrak{H}\diamond\mathfrak{N}\circ\mathfrak{F}$, $G=G/G_\mathfrak{H}\in \mathfrak{N}\circ\mathfrak{F}$. Then $G^\mathfrak{F}\in \mathfrak{N}$, and so $Z(G^\mathfrak{F})>1$. As $\mathfrak{F}=\mathtt{S}\mathfrak{F}$, we have that $H^\mathfrak{F}\leq G^\mathfrak{F}$ for every subgroup $H$ of $G$ by Lemma 2.1(2). It follows that $\mathcal{N}_{\mathfrak{H},\mathfrak{F}}(G)\geq Z(G^\mathfrak{F})>1$.\par
\medskip
\noindent\textbf{Lemma 2.4.} Let $f$ be a subgroup functor assigning to every group $G$ a characteristic subgroup $f(G)$ of $G$. Define a subgroup functor $f_i$ as follows: for every group $G$, $f_0(G)=1$; $f_i(G)/f_{i-1}(G)=f(G/f_{i-1}(G))$ for $i=1,2,\cdots$.\par
(1) If $f(G)N/N\leq f(G/N)$ for every group $G$ and every normal subgroup $N$ of $G$, then $f_i(G)N/N\leq f_i(G/N)$ for every group $G$ and every normal subgroup $N$ of $G$.\par
(2) If $f(G)N/N\leq f(G/N)$ and $f(G)\cap N\leq f(N)$ for every group $G$ and every normal subgroup $N$ of $G$, then $f_i(G)\cap N\leq f_i(N)$ for every group $G$ and every normal subgroup $N$ of $G$.\par
(3) If $f(G)N/N\leq f(G/N)$ for every group $G$ and every normal subgroup $N$ of $G$, and $f(G)\cap H\leq f(H)$ for every group $G$ and every subgroup $H$ of $G$, then $f_i(G)\cap H\leq f_i(H)$ for every group $G$ and every subgroup $H$ of $G$.\par
\medskip
\noindent\textit{Proof.} (1) By induction, we may suppose that $f_{i-1}(G)N/N\leq f_{i-1}(G/N)$. Let $f_{i-1}(G/N)=A_{i-1}/N$ and $f_{i}(G/N)=A_{i}/N$. Then $f_{i-1}(G)\leq A_{i-1}$ and $A_i/A_{i-1}=f(G/A_{i-1})$. It follows that $(f_i(G)A_{i-1}/f_{i-1}(G))/(A_{i-1}/f_{i-1}(G))\leq f((G/f_{i-1}(G))/(A_{i-1}/f_{i-1}(G)))=(A_i/f_{i-1}(G))/\linebreak[4](A_{i-1}/f_{i-1}(G))$. Therefore, $f_i(G)\leq A_i$, and so $f_i(G)N/N\leq f_i(G/N)$.\par
(2) By induction, we may assume that $f_{i-1}(G)\cap N\leq f_{i-1}(N)$. Let $f_{i-1}(G)\cap N=C_{i-1}$ and $f_i(G)\cap N=C_i$. Then $f(N/C_{i-1})(f_{i-1}(N)/C_{i-1})/(f_{i-1}(N)/C_{i-1})\leq f((N/C_{i-1})/(f_{i-1}(N)/C_{i-1}))\linebreak[4]=(f_{i}(N)/C_{i-1})/(f_{i-1}(N)/C_{i-1})$. This implies that $f(N/C_{i-1})\leq f_{i}(N)/C_{i-1}$. Clearly, $C_{i}f_{i-1}(G)/f_{i-1}(G)=f(G/f_{i-1}(G))\cap (f_{i-1}(G)N/f_{i-1}(G))\leq f(f_{i-1}(G)N/f_{i-1}(G))$. It follows that $C_{i}/C_{i-1}\leq f(N/C_{i-1})\leq f_{i}(N)/C_{i-1}$. Therefore, $C_{i}\leq f_i(N)$, and so $f_i(G)\cap N\leq f_i(N)$.\par
The proof of statement (3) is similar to (2).\par
\medskip
\noindent\textbf{Lemma 2.5.} Let $\mathfrak{H}$ be a Fitting class and $\mathfrak{F}$ a formation. Suppose that $H\leq G$ and $N\unlhd G$. Then:\par
(1) If $\mathfrak{H}=\mathtt{Q}\mathfrak{H}$, then $\mathcal{N}_{\mathfrak{H},\mathfrak{F}}^{\infty}(G)\cap N\leq \mathcal{N}_{\mathfrak{H},\mathfrak{F}}^{\infty}(N)$.\par
(2) If $\mathfrak{H}$ is a Fitting formation such that $\mathfrak{H}=\mathtt{S}\mathfrak{H}$, then $\mathcal{N}_{\mathfrak{H},\mathfrak{F}}^{\infty}(G)\cap H\leq \mathcal{N}_{\mathfrak{H},\mathfrak{F}}^{\infty}(H)$.\par
(3) If $\mathfrak{H}=\mathtt{Q}\mathfrak{H}$, then $\mathcal{N}_{\mathfrak{H},\mathfrak{F}}^{\infty}(G)N/N\leq \mathcal{N}_{\mathfrak{H},\mathfrak{F}}^{\infty}(G/N)$.\par
(4) If $\mathfrak{H}=\mathtt{Q}\mathfrak{H}$ and $N\leq \mathcal{N}_{\mathfrak{H},\mathfrak{F}}^{\infty}(G)$, then $\mathcal{N}_{\mathfrak{H},\mathfrak{F}}^{\infty}(G/N)=\mathcal{N}_{\mathfrak{H},\mathfrak{F}}^{\infty}(G)/N$.\par
(5) If $\mathfrak{H}=\mathtt{Q}\mathfrak{H}$, then $\mathcal{N}_{\mathfrak{H},\mathfrak{F}}^{\infty}(G)=\bigcap\{N|N\unlhd G, \mathcal{N}_{\mathfrak{H},\mathfrak{F}}(G/N)=1\}$.\par
\medskip
\noindent\textit{Proof.} Statements (1)-(3) directly follow from Lemmas 2.3 and 2.4.\par
(4) By definition and (3), we have that $\mathcal{N}_{\mathfrak{H},\mathfrak{F}}^{\infty}(G/N)/(\mathcal{N}_{\mathfrak{H},\mathfrak{F}}^{\infty}(G)/N)\leq \mathcal{N}_{\mathfrak{H},\mathfrak{F}}^{\infty}((G/N)/(\mathcal{N}_{\mathfrak{H},\mathfrak{F}}^{\infty}(G)/N))\linebreak[4]=1$. Therefore, $\mathcal{N}_{\mathfrak{H},\mathfrak{F}}^{\infty}(G/N)=\mathcal{N}_{\mathfrak{H},\mathfrak{F}}^{\infty}(G)/N$.\par
(5) By definition, $\mathcal{N}_{\mathfrak{H},\mathfrak{F}}(G/\mathcal{N}_{\mathfrak{H},\mathfrak{F}}^{\infty}(G))=1$. On the other hand, if $N\unlhd G$ such that $\mathcal{N}_{\mathfrak{H},\mathfrak{F}}(G/N)=1$, then $\mathcal{N}_{\mathfrak{H},\mathfrak{F}}^{\infty}(G/N)=1$. Hence by (3), $\mathcal{N}_{\mathfrak{H},\mathfrak{F}}^{\infty}(G)\leq N$. Therefore, we have that $\mathcal{N}_{\mathfrak{H},\mathfrak{F}}^{\infty}(G)=\bigcap\{N|N\unlhd G, \mathcal{N}_{\mathfrak{H},\mathfrak{F}}(G/N)=1\}$.\par
\medskip
\noindent\textbf{Lemma 2.6.} Let $\mathfrak{H}$ be a Fitting class and $\mathfrak{F}$ a formation such that $\mathfrak{F}\subseteq \mathfrak{S}$. Suppose that $G_1$ and $G_2$ are groups with $(|G_1|,|G_2|)=1$. Then $\mathcal{N}_{\mathfrak{H},\mathfrak{F}}(G_1\times G_2)=\mathcal{N}_{\mathfrak{H},\mathfrak{F}}(G_1)\times \mathcal{N}_{\mathfrak{H},\mathfrak{F}}(G_2)$ and $\mathcal{N}_{\mathfrak{H},\mathfrak{F}}^{\infty}(G_1\times G_2)=\mathcal{N}_{\mathfrak{H},\mathfrak{F}}^{\infty}(G_1)\times \mathcal{N}_{\mathfrak{H},\mathfrak{F}}^{\infty}(G_2)$.\par
\medskip
\noindent\textit{Proof.} We only need to prove that $\mathcal{N}_{\mathfrak{H},\mathfrak{F}}(G_1\times G_2)=\mathcal{N}_{\mathfrak{H},\mathfrak{F}}(G_1)\times \mathcal{N}_{\mathfrak{H},\mathfrak{F}}(G_2)$. Let $G=G_1\times G_2$. Since $(|G_1|,|G_2|)=1$, for every subgroup $H$ of $G$, we have that $H=(H\cap G_1)\times (H\cap G_2)$. By \cite[Chap. IV, Theorem 1.18]{Doe}, $H^\mathfrak{F}=(H\cap G_1)^\mathfrak{F}\times (H\cap G_2)^\mathfrak{F}$. Then it is easy to see that $G_\mathfrak{H}={(G_1)}_\mathfrak{H}\times {(G_2)}_\mathfrak{H}$, and so $H^\mathfrak{F}G_\mathfrak{H}=(H\cap G_1)^\mathfrak{F}{(G_1)}_\mathfrak{H}\times (H\cap G_2)^\mathfrak{F}{(G_2)}_\mathfrak{H}$. This implies that $N_G(H^\mathfrak{F}G_\mathfrak{H})=N_{G_1}((H\cap G_1)^\mathfrak{F}{(G_1)}_\mathfrak{H})\times N_{G_2}((H\cap G_2)^\mathfrak{F}{(G_2)}_\mathfrak{H})$. Hence $\mathcal{N}_{\mathfrak{H},\mathfrak{F}}(G)=\bigcap_{H\leq G}N_G(H^{\mathfrak{F}}G_{\mathfrak{H}})=\bigcap_{H\leq G}N_{G_1}((H\cap G_1)^{\mathfrak{F}}{(G_1)}_{\mathfrak{H}})\times \bigcap_{H\leq G}N_{G_2}((H\cap G_2)^{\mathfrak{F}}{(G_2)}_{\mathfrak{H}})=\mathcal{N}_{\mathfrak{H},\mathfrak{F}}(G_1)\times \mathcal{N}_{\mathfrak{H},\mathfrak{F}}(G_2)$.\par
\medskip
\noindent\textbf{Lemma 2.7.} Let $\mathfrak{F}$ be a formation. Then $\mathfrak{F}$ satisfies the $\pi$-boundary condition (I) if and only if $\mbox{Crit}_\mathtt{S}(\mathfrak{F})\subseteq \mathfrak{S}_\pi\circ\mathfrak{F}$.\par
\medskip
\noindent\textit{Proof.} The necessity is evident. So we only need to prove the sufficiency. Suppose that $\mbox{Crit}_\mathtt{S}(\mathfrak{F})\subseteq \mathfrak{S}_\pi\circ\mathfrak{F}$. Let $G\in \mbox{Crit}_\mathtt{S}(\mathfrak{F})$. If $G^\mathfrak{F}\leq \Phi(G)$, then there is nothing to prove. We may, therefore, assume that $G^\mathfrak{F}\nleq \Phi(G)$. Let $G^\mathfrak{F}/L$ be a $G$-chief factor. Clearly, $G^\mathfrak{F}/L\in \mathfrak{N}_\pi$. If $L\nleq \Phi(G)$, then $G$ has a maximal subgroup $M$ such that $G=LM$. Since $M\in \mathfrak{F}$, $G/L\cong M/L\cap M\in \mathfrak{F}$, and so $G^\mathfrak{F}\leq L$, which is absurd. Hence $L\leq \Phi(G)$. This implies that $L=G^\mathfrak{F}\cap \Phi(G)$. Since $G^\mathfrak{F}\Phi(G)/\Phi(G)\cong G^\mathfrak{F}/G^\mathfrak{F}\cap \Phi(G)\in \mathfrak{N}_\pi$, we have that $G^\mathfrak{F}\in \mathfrak{N}_\pi$ by \cite[Lemma 3.1]{Bal1}. This shows that $G\in \mathfrak{N}_\pi\circ\mathfrak{F}$, and thus $\mbox{Crit}_\mathtt{S}(\mathfrak{F})\subseteq \mathfrak{N}_\pi\circ\mathfrak{F}$.\par
\medskip
\noindent\textbf{Lemma 2.8.} Let $\mathfrak{F}$ be a saturated formation and $\pi\subseteq \pi(\mathfrak{F})$. Suppose that $H\leq G$ and $N\unlhd G$. Then:\par
(1) If $N\leq Z_{\pi\mathfrak{F}}(G)$, then $Z_{\pi\mathfrak{F}}(G/N)=Z_{\pi\mathfrak{F}}(G)/N$.\par
(2) $Z_{\pi\mathfrak{F}}(G)N/N\leq Z_{\pi\mathfrak{F}}(G/N)$.\par
(3) If $\mathfrak{F}=\mathtt{S}\mathfrak{F}$ (resp. $\mathfrak{F}=\mathtt{S_n}\mathfrak{F}$), then $Z_{\pi\mathfrak{F}}(G)\cap H\leq Z_{\pi\mathfrak{F}}(H)$ (resp. $Z_{\pi\mathfrak{F}}(G)\cap N\leq Z_{\pi\mathfrak{F}}(N)$).\par
(4) If $\mathfrak{G}_{\pi'}\circ\mathfrak{F}=\mathfrak{F}$ and $G/Z_{\pi\mathfrak{F}}(G)\in \mathfrak{F}$, then $G\in \mathfrak{F}$.\par
(5) If $\mathfrak{F}=\mathtt{S}\mathfrak{F}$ (resp. $\mathfrak{F}=\mathtt{S_n}\mathfrak{F}$), $\mathfrak{G}_{\pi'}\circ\mathfrak{F}=\mathfrak{F}$ and $H\in \mathfrak{F}$ (resp. $N\in \mathfrak{F}$), then $HZ_{\pi\mathfrak{F}}(G)\in \mathfrak{F}$ (resp. $NZ_{\pi\mathfrak{F}}(G)\in \mathfrak{F}$).\par
(6) $Z_{\pi\mathfrak{F}}(G)=Z_{\pi(\mathfrak{G}_{\pi'}\circ\mathfrak{F})}(G)$.\par
(7) If $\mathfrak{F}=\mathtt{S_n}\mathfrak{F}$, then $Z_{\pi\mathfrak{F}}(G)\in \mathfrak{G}_{\pi'}\circ\mathfrak{F}$.\par
\medskip
\noindent\textit{Proof.} Statement (1) is evident by definition.\par
Statements (2)-(5) were proved in \cite[Lemma 2.2]{Guo}.\par
(6) Let $\mathfrak{F}=LF(F)$, where $F$ is the canonical local definition of $\mathfrak{F}$. Then by \cite[Chap. IV, Theorem 3.13]{Doe}, $\mathfrak{G}_{\pi'}\circ\mathfrak{F}=LF(H)$, where $H(p)=F(p)$ for all $p\in \pi$ and $H(p)=\mathfrak{G}_{\pi'}\circ\mathfrak{F}$ for all $p\in \pi'$. Then by definition, it is easy to see that $Z_{\pi\mathfrak{F}}(G)=Z_{\pi(\mathfrak{G}_{\pi'}\circ\mathfrak{F})}(G)$.\par
Statement (7) follows from (5) and (6).\par
\medskip
\noindent\textbf{Remark 2.9.} Note that there exist several minor mistakes in \cite{Guo}. In \cite[Lemmas 2.2(6) and 2.2(7)]{Guo} and \cite[Lemma 2.4(g)]{Guo}, ``$\mathfrak{G}_\sigma\circ\mathfrak{F}=\mathfrak{F}$" should be corrected as ``$\mathfrak{G}_{\pi'}\circ\mathfrak{F}=\mathfrak{F}$"; and in \cite[Lemma 2.2(5)]{Guo}, ``$Z_{\pi\mathfrak{F}}(H)\cap A$" should be corrected as ``$Z_{\pi\mathfrak{F}}(A)\cap H$".\par
\medskip
\noindent\textbf{Lemma 2.10.}\cite[Lemma 2.5]{Ski} Let $\mathfrak{F}=LF(F)$ be a saturated formation, where $F$ is the canonical local definition of $\mathfrak{F}$, and $E$ a normal $p$-subgroup of $G$. If $E\leq Z_\mathfrak{F}(G)$, then $G/C_G(E)\in F(p)$.\par
\medskip
\noindent\textbf{Lemma 2.11.} Let $\mathfrak{F}$ be a formation and $\mathfrak{B}=\mathfrak{N}_\pi\circ \mathfrak{F}$. Then:\par
(1) $\mathfrak{B}=LF(b)$ with $b(p)=\mathfrak{F}$ for all $p\in \pi$ and $b(p)=\mathfrak{B}=\mathfrak{N}_\pi\circ \mathfrak{F}$ for all $p\in \pi'$.\par
(2) The canonical local definition $B$ of $\mathfrak{B}$ can be defined as follows: $B(p)=\mathfrak{G}_p\circ\mathfrak{F}$ for all $p\in \pi$ and $B(p)=\mathfrak{B}=\mathfrak{N}_\pi\circ\mathfrak{F}$ for all $p\in \pi'$.\par
\medskip
\noindent\textit{Proof.} Statement (1) directly follows from \cite[Lemma 1]{Ski2}, and Statement (2) follows from \cite[Chap. IV, Lemma 3.13]{Doe}.\par
\medskip
\noindent\textbf{Lemma 2.12.} Let $\mathfrak{F}$ be a formation. Then:\par
(1) $Z_{\pi(\mathfrak{N}\circ\mathfrak{F})}(G)=1$ if and only if $C_G(G^\mathfrak{F})=1$ and $O_{\pi'}(G)=1$.\par
(2) $Z_{\pi(\mathfrak{N}\circ\mathfrak{F})}(G)\cap G^\mathfrak{F}=Z_{\pi\mathfrak{N}}(G^\mathfrak{F})$.\par
(3) If $\mathfrak{F}$ is saturated, then $Z_{\pi(\mathfrak{N}\circ\mathfrak{F})}(G)/Z_{\pi\mathfrak{N}}(G^\mathfrak{F})=Z_{\pi\mathfrak{F}}(G/Z_{\pi\mathfrak{N}}(G^\mathfrak{F}))$.\par
\medskip
\noindent\textit{Proof.} (1) Suppose that $C_G(G^\mathfrak{F})=1$ and $O_{\pi'}(G)=1$. If $Z_{\pi(\mathfrak{N}\circ\mathfrak{F})}(G)>1$, then let $N$ be a minimal normal subgroup of $G$ contained in $Z_{\pi(\mathfrak{N}\circ\mathfrak{F})}(G)$. Clearly, $N$ is not a $\pi'$-group. Then by Lemma 2.11(1), we have that $G/C_G(N)\in \mathfrak{F}$, and so $N\leq C_G(G^\mathfrak{F})=1$, a contradiction. Thus $Z_{\pi(\mathfrak{N}\circ\mathfrak{F})}(G)=1$. Now assume that $Z_{\pi(\mathfrak{N}\circ\mathfrak{F})}(G)=1$. Then clearly, $O_{\pi'}(G)=1$. Suppose that $C_G(G^\mathfrak{F})>1$, and let $N$ be a minimal normal subgroup of $G$ contained in $C_G(G^\mathfrak{F})$. Then $G^\mathfrak{F}\leq C_G(N)$ and $N$ is not a $\pi'$-group. Hence by Lemma 2.11(1) again, $N\leq Z_{\pi(\mathfrak{N}\circ\mathfrak{F})}(G)$, which is impossible. Therefore, $C_G(G^\mathfrak{F})=1$.\par
(2) Firstly, we prove that $Z_{\pi\mathfrak{N}}(G^\mathfrak{F})\leq Z_{\pi(\mathfrak{N}\circ\mathfrak{F})}(G)$. If $Z_{\pi(\mathfrak{N}\circ\mathfrak{F})}(G)>1$, then by induction, $Z_{\pi\mathfrak{N}}((G/Z_{\pi(\mathfrak{N}\circ\mathfrak{F})}(G))^\mathfrak{F})\leq Z_{\pi(\mathfrak{N}\circ\mathfrak{F})}(G/Z_{\pi(\mathfrak{N}\circ\mathfrak{F})}(G))=1$. By Lemmas 2.1(1) and 2.8(2), $Z_{\pi\mathfrak{N}}(G^\mathfrak{F})\leq Z_{\pi(\mathfrak{N}\circ\mathfrak{F})}(G)$. We may, therefore, assume that $Z_{\pi(\mathfrak{N}\circ\mathfrak{F})}(G)=1$. Then by (1), $C_G(G^\mathfrak{F})=1$ and $O_{\pi'}(G)=1$. It follows that $Z(G^\mathfrak{F})=1$ and $O_{\pi'}(G^\mathfrak{F})=1$. By (1) again, $Z_{\pi\mathfrak{N}}(G^\mathfrak{F})=1$. Consequently, $Z_{\pi\mathfrak{N}}(G^\mathfrak{F})\leq Z_{\pi(\mathfrak{N}\circ\mathfrak{F})}(G)$.\par
Suppose that $Z_{\pi\mathfrak{N}}(G^\mathfrak{F})>1$. Then by induction and Lemma 2.1(1), $Z_{\pi(\mathfrak{N}\circ\mathfrak{F})}(G/Z_{\pi\mathfrak{N}}(G^\mathfrak{F}))\cap (G^\mathfrak{F}/Z_{\pi\mathfrak{N}}(G^\mathfrak{F}))=Z_{\pi\mathfrak{N}}(G^\mathfrak{F}/Z_{\pi\mathfrak{N}}(G^\mathfrak{F}))=1$. Hence by Lemma 2.8(1), $Z_{\pi(\mathfrak{N}\circ\mathfrak{F})}(G)\cap G^\mathfrak{F}=Z_{\pi\mathfrak{N}}(G^\mathfrak{F})$. We may, therefore, assume that $Z_{\pi\mathfrak{N}}(G^\mathfrak{F})=1$. Then by (1), $Z(G^\mathfrak{F})=1$ and $O_{\pi'}(G^\mathfrak{F})=1$. If $Z_{\pi(\mathfrak{N}\circ\mathfrak{F})}(G)\cap G^\mathfrak{F}>1$, then let $N$ be a minimal normal subgroup of $G$ contained in $Z_{\pi(\mathfrak{N}\circ\mathfrak{F})}(G)\cap G^\mathfrak{F}$. Since $O_{\pi'}(G^\mathfrak{F})=1$, $N$ is not a $\pi'$-group. It follows from Lemma 2.11(1) that $G/C_G(N)\in \mathfrak{F}$, and so $N\leq Z(G^\mathfrak{F})$, a contradiction. Therefore, $Z_{\pi(\mathfrak{N}\circ\mathfrak{F})}(G)\cap G^\mathfrak{F}=1$.\par
(3) If $Z_{\pi\mathfrak{N}}(G^\mathfrak{F})>1$, then by induction, $Z_{\pi(\mathfrak{N}\circ\mathfrak{F})}(G/Z_{\pi\mathfrak{N}}(G^\mathfrak{F}))=Z_{\pi\mathfrak{F}}(G/Z_{\pi\mathfrak{N}}(G^\mathfrak{F}))$. Hence by (2) and Lemma 2.8(1), $Z_{\pi(\mathfrak{N}\circ\mathfrak{F})}(G)/Z_{\pi\mathfrak{N}}(G^\mathfrak{F})=Z_{\pi\mathfrak{F}}(G/Z_{\pi\mathfrak{N}}(G^\mathfrak{F}))$. We may, therefore, assume that $Z_{\pi\mathfrak{N}}(G^\mathfrak{F})=1$. Then by (2), $Z_{\pi(\mathfrak{N}\circ\mathfrak{F})}(G)\cap G^\mathfrak{F}=1$, and so $Z_{\pi(\mathfrak{N}\circ\mathfrak{F})}(G)\leq C_G(G^\mathfrak{F})$. By \cite[Chap. IV, Theorem 6.13]{Doe}, $C_G(G^\mathfrak{F})=Z_\mathfrak{F}(G)\leq Z_{\pi\mathfrak{F}}(G)\leq Z_{\pi(\mathfrak{N}\circ\mathfrak{F})}(G)$. This implies that $Z_{\pi(\mathfrak{N}\circ\mathfrak{F})}(G)=Z_{\pi\mathfrak{F}}(G)$.\par
\medskip
\noindent\textbf{Lemma 2.13.}\cite[Lemma 2.10]{Ski} Let $\mathfrak{F}=LF(F)$ be a saturated formation with $p\in \pi(\mathfrak{F})$, where $F$ is the canonical local definition of $\mathfrak{F}$. Suppose that $G$ is a group of minimal order in the set of all groups $G\in \mbox{Crit}_\mathtt{S}(F(p))$ and $G\notin \mathfrak{F}$. Then $G^\mathfrak{F}$ is the unique minimal normal subgroup of $G$ and $O_p(G)=\Phi(G)=1$.\par
\section{Proofs of Main Results}
\noindent\textbf{Lemma 3.1.} Let $\mathfrak{H}$ be a saturated Fitting formation such that $\mathfrak{G}_{\pi'}\subseteq\mathfrak{H}=\mathtt{E}\mathfrak{H}$ and $\mathfrak{F}$ a formation. Then $\mathcal{N}_{\mathfrak{H},\mathfrak{F}}^{\infty}(G)\in \mathfrak{H}\diamond\mathfrak{N}\circ\mathfrak{F}$ if one of the following holds:\par
(i) $\mathfrak{F}=\mathtt{S_n}\mathfrak{F}$ and $G^{\mathfrak{H}\diamond\mathfrak{N}\circ\mathfrak{F}}\in \mathfrak{S_{\pi}}$.\par
(ii) $\mathfrak{F}$ satisfies the $\pi$-boundary condition (I).\par
\medskip
\noindent\textit{Proof.} Assume that the result is false and let $G$ be a counterexample of minimal order. Note that if the condition (i) holds, then since $\mathfrak{H}\diamond\mathfrak{N}\circ\mathfrak{F}=\mathtt{S_n}(\mathfrak{H}\diamond\mathfrak{N}\circ\mathfrak{F})$, we have that $(\mathcal{N}_{\mathfrak{H},\mathfrak{F}}^{\infty}(G))^{\mathfrak{H}\diamond\mathfrak{N}\circ\mathfrak{F}}\leq G^{\mathfrak{H}\diamond\mathfrak{N}\circ\mathfrak{F}}\in \mathfrak{S}_\pi$ by Lemma 2.1(2). Hence the condition (i) holds for $\mathcal{N}_{\mathfrak{H},\mathfrak{F}}^{\infty}(G)$ when the condition (i) holds for $G$. If $\mathcal{N}_{\mathfrak{H},\mathfrak{F}}^{\infty}(G)<G$, then by the choice of $G$ and Lemma 2.5(1), $\mathcal{N}_{\mathfrak{H},\mathfrak{F}}^{\infty}(G)=\mathcal{N}_{\mathfrak{H},\mathfrak{F}}^{\infty}(\mathcal{N}_{\mathfrak{H},\mathfrak{F}}^{\infty}(G))\in \mathfrak{H}\diamond\mathfrak{N}\circ\mathfrak{F}$, a contradiction. We may, therefore, assume that $\mathcal{N}_{\mathfrak{H},\mathfrak{F}}^{\infty}(G)=G$. Let $N$ be any minimal normal subgroup of $G$. Then by Lemma 2.1(1), the condition (i) holds for $G/N$ when the condition (i) holds for $G$. Hence by the choice of $G$ and Lemma 2.5(3), $G/N=\mathcal{N}_{\mathfrak{H},\mathfrak{F}}^{\infty}(G/N)\in \mathfrak{H}\diamond\mathfrak{N}\circ\mathfrak{F}$. Clearly, $\mathfrak{H}\diamond\mathfrak{N}\circ\mathfrak{F}$ is a saturated formation by \cite[Chap. IV, Theorem 4.8]{Doe}. This implies that $N$ is the unique minimal normal subgroup of $G$.\par
If $G_\mathfrak{H}>1$, then $N\leq G_\mathfrak{H}$. By Lemmas 2.2(3) and 2.2(4), $(G/N)_\mathfrak{H}=G_\mathfrak{H}/N$. Since $G/N\in \mathfrak{H}\diamond\mathfrak{N}\circ\mathfrak{F}$, $G/G_\mathfrak{H}\in \mathfrak{N}\circ\mathfrak{F}$, and thus $G\in \mathfrak{H}\diamond\mathfrak{N}\circ\mathfrak{F}$, a contradiction. Therefore, $G_\mathfrak{H}=1$, and so $O_{\pi'}(G)=1$. If $N\leq \Phi(G)$, then $G\in \mathfrak{H}\diamond\mathfrak{N}\circ\mathfrak{F}$, which is impossible. Hence $N\nleq \Phi(G)$. It follows that $G$ has a maximal subgroup $M$ such that $N\nleq M$. Since $\mathcal{N}_{\mathfrak{H},\mathfrak{F}}(G)>1$ and $N$ is the unique minimal normal subgroup of $G$, we have that $N\leq \mathcal{N}_{\mathfrak{H},\mathfrak{F}}(G)$. Then by the definition of $\mathcal{N}_{\mathfrak{H},\mathfrak{F}}(G)$, $N\leq N_G(M^\mathfrak{F})$. This induces that $M^\mathfrak{F}\unlhd G$. Hence $M^\mathfrak{F}=1$, and so $M\in \mathfrak{F}$. It follows that $G/N\cong M/N\cap M\in \mathfrak{F}$, and thereby $G^\mathfrak{F}\leq N$. Since $1<G^{\mathfrak{H}\diamond\mathfrak{N}\circ\mathfrak{F}}\leq G^\mathfrak{F}\leq N$, $N=G^{\mathfrak{H}\diamond\mathfrak{N}\circ\mathfrak{F}}=G^\mathfrak{F}$.\par
We claim that $N\in \mathfrak{N}$. If the condition (i) holds, then $N\in \mathfrak{S_{\pi}}$. As $O_{\pi'}(G)=1$, $N\in \mathfrak{N}$. Now assume that the condition (ii) holds. Then since $G\notin \mathfrak{F}$, we may take a subgroup $K$ of $G$ such that $K\in \mbox{Crit}_\mathtt{S}(\mathfrak{F})\subseteq \mathfrak{N_{\pi}}\circ\mathfrak{F}$. If $N\notin \mathfrak{N}$, then $C_G(N)=1$. Since $N\leq \mathcal{N}_{\mathfrak{H},\mathfrak{F}}(G)$ and $G_\mathfrak{H}=1$, we have that $N\leq N_G(K^\mathfrak{F})$, and so $N\cap K^\mathfrak{F}\unlhd N$. As $K^\mathfrak{F}\in \mathfrak{N_{\pi}}$ and $O_{\pi'}(G)=1$, $K^\mathfrak{F}\in \mathfrak{N}$. By \cite[Chap. A, Proposition 4.13(b)]{Doe}, $N\cap K^\mathfrak{F}=1$. It follows that $K^\mathfrak{F}\leq C_G(N)=1$, and thus $K\in \mathfrak{F}$, a contradiction. Hence $N\in \mathfrak{N}$. Therefore, our claim holds. This induces that $G\in \mathfrak{N}\circ\mathfrak{F}\subseteq \mathfrak{H}\diamond\mathfrak{N}\circ\mathfrak{F}$. The final contradiction completes the proof.\par
\medskip
\noindent\textbf{Proof of Theorem A.} It is obvious that (1) implies (2) and (2) implies (3). Suppose that (3) holds, that is, $G/\mathcal{N}_{\mathfrak{H},\mathfrak{F}}^{\infty}(G)\in \mathfrak{H}\diamond\mathfrak{N}\circ\mathfrak{F}$. If $\mathcal{N}_{\mathfrak{H},\mathfrak{F}}(G/N)=1$ for some proper normal subgroup $N$ of $G$, then by Lemma 2.5(5), $\mathcal{N}_{\mathfrak{H},\mathfrak{F}}^{\infty}(G)\leq N$, and so $G/N\in \mathfrak{H}\diamond\mathfrak{N}\circ\mathfrak{F}$. Hence by Lemma 2.3(4), either $G=N$ or $\mathcal{N}_{\mathfrak{H},\mathfrak{F}}(G/N)>1$, a contradiction. This induces that (3) implies (4). Now assume that (4) holds. Then since $\mathcal{N}_{\mathfrak{H},\mathfrak{F}}(G/\mathcal{N}_{\mathfrak{H},\mathfrak{F}}^{\infty}(G))=1$, we have that $G=\mathcal{N}_{\mathfrak{H},\mathfrak{F}}^{\infty}(G)$. Hence (4) implies (5). Finally, by Lemma 3.1, we get that (5) implies (1). This finishes the proof of the theorem.\par
\medskip
Since $\mathfrak{N}_\pi=\mathfrak{G}_{\pi'}\circ\mathfrak{N}=\mathfrak{G}_{\pi'}\diamond\mathfrak{N}$, the next corollary directly follows from Theorem A, which is also a generalization of \cite[Theorem A]{Su} and \cite[Theorem B]{Su}.\par
\medskip
\noindent\textbf{Corollary 3.2.} Let $\mathfrak{F}$ be a formation such that $\mathfrak{F}=\mathtt{S}\mathfrak{F}$. Suppose that one of the following holds:\par
(i) $G\in \mathfrak{S_{\pi}}\circ \mathfrak{F}$.\par
(ii) $\mathfrak{F}$ satisfies the $\pi$-boundary condition (I).\par
\noindent Then the following statements are equivalent:\par
(1) $G\in \mathfrak{N}_\pi\circ\mathfrak{F}$.\par
(2) $G/\mathcal{N}_{\pi'\mathfrak{F}}(G)\in \mathfrak{N}_\pi\circ\mathfrak{F}$.\par
(3) $G/\mathcal{N}_{\pi'\mathfrak{F}}^{\infty}(G)\in \mathfrak{N}_\pi\circ\mathfrak{F}$.\par
(4) $\mathcal{N}_{\pi'\mathfrak{F}}(G/N)>1$ for every proper normal subgroup $N$ of $G$.\par
(5) $G=\mathcal{N}_{\pi'\mathfrak{F}}^{\infty}(G)$.\par
\medskip
In the sequel of this section, we restrict our attention to $\pi\mathfrak{F}$-norms.\par
\medskip
\noindent\textbf{Lemma 3.3.} Let $\mathfrak{F}$ be a formation such that $\mathfrak{F}=\mathtt{S}\mathfrak{F}$. Then $Z_{\pi(\mathfrak{N}\circ\mathfrak{F})}(G)\leq \mathcal{N}_{\pi'\mathfrak{F}}^{\infty}(G)$.\par
\medskip
\noindent\textit{Proof.} If $\mathcal{N}_{\pi'\mathfrak{F}}(G)>1$, then by induction, $Z_{\pi(\mathfrak{N}\circ\mathfrak{F})}(G/\mathcal{N}_{\pi'\mathfrak{F}}^{\infty}(G))\leq \mathcal{N}_{\pi'\mathfrak{F}}^{\infty}(G/\mathcal{N}_{\pi'\mathfrak{F}}^{\infty}(G))=1$. By Lemma 2.8(2), $Z_{\pi(\mathfrak{N}\circ\mathfrak{F})}(G)\mathcal{N}_{\pi'\mathfrak{F}}^{\infty}(G)/\mathcal{N}_{\pi'\mathfrak{F}}^{\infty}(G)=1$, and so $Z_{\pi(\mathfrak{N}\circ\mathfrak{F})}(G)\leq \mathcal{N}_{\pi'\mathfrak{F}}^{\infty}(G)$. Hence we may assume that $\mathcal{N}_{\pi'\mathfrak{F}}(G)=1$. Since $\mathfrak{F}=\mathtt{S}\mathfrak{F}$, $H^\mathfrak{F}\leq G^\mathfrak{F}$ for every subgroup $H$ of $G$ by Lemma 2.1(2), and thereby $C_G(G^\mathfrak{F})O_{\pi'}(G)\leq N_G(H^\mathfrak{F}O_{\pi'}(G))$. This implies that $C_G(G^\mathfrak{F})=1$ and $O_{\pi'}(G)=1$. Then by Lemma 2.12(1), $1=Z_{\pi(\mathfrak{N}\circ\mathfrak{F})}(G)\leq \mathcal{N}_{\pi'\mathfrak{F}}^{\infty}(G)$.\par
\medskip
\noindent\textbf{Proofs of Theorem B(1) and Theorem D.} We need to prove that if either $G\in \mathfrak{S_{\pi}}\circ \mathfrak{F}$ or $\mathfrak{F}$ satisfies the $\pi$-boundary condition (II), then $\mathcal{N}_{\pi'\mathfrak{F}}^{\infty}(G)=Z_{\pi(\mathfrak{N}\circ\mathfrak{F})}(G)$. Suppose that the result is false and let $L$ be a counterexample of minimal order. By Lemma 3.3, $Z_{\pi(\mathfrak{N}\circ\mathfrak{F})}(L)\leq \mathcal{N}_{\pi'\mathfrak{F}}^{\infty}(L)$. We may, therefore, assume that $\mathcal{N}_{\pi'\mathfrak{F}}(L)>1$. If $Z_{\pi(\mathfrak{N}\circ\mathfrak{F})}(L)>1$, then by the choice of $L$ and Lemma 2.8(1), $\mathcal{N}_{\pi'\mathfrak{F}}^{\infty}(L/Z_{\pi(\mathfrak{N}\circ\mathfrak{F})}(L))=Z_{\pi(\mathfrak{N}\circ\mathfrak{F})}(L/Z_{\pi(\mathfrak{N}\circ\mathfrak{F})}(L))=1$. Hence by Lemma 2.5(4), $\mathcal{N}_{\pi'\mathfrak{F}}^{\infty}(L)=Z_{\pi(\mathfrak{N}\circ\mathfrak{F})}(L)$, a contradiction. Therefore, $Z_{\pi(\mathfrak{N}\circ\mathfrak{F})}(L)=1$, and thereby $O_{\pi'}(L)=1$.\par
Now let $N$ be any minimal normal subgroup of $L$ contained in $\mathcal{N}_{\pi'\mathfrak{F}}(L)$. If $L$ has a minimal normal subgroup $R$ which is different from $N$, then by the choice of $L$, $\mathcal{N}_{\pi'\mathfrak{F}}^{\infty}(L/R)=Z_{\pi(\mathfrak{N}\circ\mathfrak{F})}(L/R)$. It follows from Lemma 2.5(3) that $NR/R\leq \mathcal{N}_{\pi'\mathfrak{F}}^{\infty}(L/R)=Z_{\pi(\mathfrak{N}\circ\mathfrak{F})}(L/R)$. By $L$-isomorphism $N\cong NR/R$, we have that $N\leq Z_{\pi(\mathfrak{N}\circ\mathfrak{F})}(L)$, which is absurd. Thus $N$ is the unique minimal normal subgroup of $L$. If $N\nleq \Phi(L)$, then $L$ has a maximal subgroup $M$ such that $N\nleq M$. Since $N\leq \mathcal{N}_{\pi'\mathfrak{F}}(L)$, $N\leq N_L(M^\mathfrak{F})$ for $O_{\pi'}(L)=1$. This implies that $M^\mathfrak{F}\unlhd L$, and so $M^\mathfrak{F}=1$. Therefore, $M\in \mathfrak{F}$. Then $L/N\cong M/N\cap M\in \mathfrak{F}$. As $N\leq \mathcal{N}_{\pi'\mathfrak{F}}(L)$, $L/\mathcal{N}_{\pi'\mathfrak{F}}(L)\in \mathfrak{F}$. Note that $\mathfrak{F}$ satisfies the $\pi$-boundary condition (I) if $\mathfrak{F}$ satisfies the $\pi$-boundary condition (II). By Corollary 3.2, we have that $L\in \mathfrak{N}_{\pi}\circ \mathfrak{F}$. This induces that $L=Z_{\pi(\mathfrak{N}_\pi\circ\mathfrak{F})}(L)=Z_{\pi(\mathfrak{N}\circ\mathfrak{F})}(L)$ by Lemma 2.8(6), a contradiction. Hence $N\leq \Phi(L)$, and so $N$ is an elementary abelian $p$-group with $p\in \pi$. Let $M$ be any maximal subgroup of $L$. By the choice of $L$, $\mathcal{N}_{\pi'\mathfrak{F}}^{\infty}(M)=Z_{\pi(\mathfrak{N}\circ\mathfrak{F})}(M)$. Then by Lemma 2.5(2), $N\leq \mathcal{N}_{\pi'\mathfrak{F}}^{\infty}(L)\cap M\leq \mathcal{N}_{\pi'\mathfrak{F}}^{\infty}(M)=Z_{\pi(\mathfrak{N}\circ\mathfrak{F})}(M)$. Thus $N\leq Z_{\mathfrak{N}\circ\mathfrak{F}}(M)$. By Lemmas 2.10 and 2.11(2), $M/C_M(N)\in \mathfrak{G}_p\circ\mathfrak{F}$. If $C_L(N)\nleq M$, then $L=C_L(N)M$, and so $L/C_L(N)\cong M/C_M(N)\in \mathfrak{G}_p\circ\mathfrak{F}$. This shows that $N\leq Z_{\pi(\mathfrak{N}\circ\mathfrak{F})}(L)$ by Lemma 2.11(2), which is impossible. Hence $C_L(N)\leq M$, and thereby $C_L(N)\leq \Phi(L)$. Since $N$ is the unique minimal normal subgroup of $L$, $\Phi(L)$ is a $p$-subgroup of $L$. Therefore, $C_L(N)$ is also a $p$-subgroup of $L$. This implies that $M\in \mathfrak{G}_p\circ\mathfrak{F}$. If $L\in \mathfrak{G}_p\circ\mathfrak{F}$, then $L\in \mathfrak{N}_{\pi}\circ \mathfrak{F}$, a contradiction as above. Hence $L\in \mbox{Crit}_\mathtt{S}(\mathfrak{G}_p\circ\mathfrak{F})$. Then, in both cases, $L\in \mathfrak{S_{\pi}}\circ \mathfrak{F}$.\par
Let $F_p(L)$ be the $p$-Fitting subgroup of $L$, that is, the $\mathfrak{N}_p$-radical of $L$. As $N$ is the unique minimal normal subgroup of $L$, we have that $O_{p'}(L)=1$, and so $F_p(L)=O_p(L)$. By \cite[Chap. A, Theorem 13.8(a)]{Doe}, $F_p(L)\leq C_L(N)\leq \Phi(L)$. This induces that $F_p(L)=\Phi(L)$. Since $L\in \mathfrak{S_{\pi}}\circ \mathfrak{F}$, $L^\mathfrak{F}\in \mathfrak{S_{\pi}}$. If $L^\mathfrak{F}\leq \Phi(L)$, then $L\in \mathfrak{G}_p\circ\mathfrak{F}$, which is absurd. Thus $L^\mathfrak{F}\nleq \Phi(L)$. Let $A/\Phi(L)$ be an $L$-chief factor contained in $L^\mathfrak{F}\Phi(L)/\Phi(L)$. Then $A/\Phi(L)\in \mathfrak{N_{\pi}}$, and so $A/\Phi(L)\in \mathfrak{N}_p$. Hence by \cite[Lemma 3.1]{Bal1}, we have that $A\in \mathfrak{N}_p$. This implies that $A\leq F_p(L)=\Phi(L)$, a contradiction. The proof is thus completed.\par
\medskip
\noindent\textbf{Proofs of Theorems B(2) and B(3).} (2) Suppose that $\mathcal{N}_{\pi'\mathfrak{F}}^{\infty}(G)=Z_{\pi(\mathfrak{N}\circ\mathfrak{F})}(G)$ holds for every group $G$. Obviously, for any group $G\in \mbox{Crit}_\mathtt{S}(\mathfrak{F})$, $G=\mathcal{N}_{\pi'\mathfrak{F}}(G)=\mathcal{N}_{\pi'\mathfrak{F}}^{\infty}(G)$. It follows that $G=Z_{\pi(\mathfrak{N}\circ\mathfrak{F})}(G)$, and so $G\in \mathfrak{N}_\pi\circ\mathfrak{F}$ by Lemma 2.8(7). Hence $\mathfrak{F}$ satisfies the $\pi$-boundary condition (I).\par
(3) The necessity is obvious. So we only need to prove the sufficiency. For any group $G\in \mbox{Crit}_\mathtt{S}(\mathfrak{F})$ and any $p\in \pi$, either $G\in \mathfrak{S}_\pi\circ \mathfrak{F}$ or $G\in \mbox{Crit}_\mathtt{S}(\mathfrak{G}_p\circ\mathfrak{F})\backslash (\mathfrak{S}_\pi\circ \mathfrak{F})$. In the former case, $G\in \mathfrak{N}_\pi\circ \mathfrak{F}$ by Lemma 2.7. In the latter case, a same discussion as in the proof of (2) shows that $G\in \mathfrak{N}_\pi\circ\mathfrak{F}$. Therefore, $\mathfrak{F}$ satisfies the $\pi$-boundary condition (I). The rest of the proof is similar to the proofs of Theorem B(1) and Theorem D.\par
\medskip
\noindent\textbf{Corollary 3.4.} Let $\mathfrak{F}$ be a formation such that $\mathfrak{F}=\mathtt{S}\mathfrak{F}$. Suppose that one of the following holds:\par
(i) $G\in \mathfrak{S}_\pi\circ\mathfrak{F}$.\par
(ii) $\mathfrak{F}$ satisfies the $\pi$-boundary condition (II).\par
\noindent Then $\mathcal{N}_{\pi'\mathfrak{F}}^{\infty}(G)/Z_{\pi\mathfrak{N}}(G^\mathfrak{F})=\mathcal{N}_{\pi'\mathfrak{F}}(G/Z_{\pi\mathfrak{N}}(G^\mathfrak{F}))=Z_{\pi(\mathfrak{N}\circ\mathfrak{F})}(G)/Z_{\pi\mathfrak{N}}(G^\mathfrak{F})$. In particular, if $\mathfrak{F}$ is saturated, then $\mathcal{N}_{\pi'\mathfrak{F}}^{\infty}(G)/Z_{\pi\mathfrak{N}}(G^\mathfrak{F})=\mathcal{N}_{\pi'\mathfrak{F}}(G/Z_{\pi\mathfrak{N}}(G^\mathfrak{F}))=Z_{\pi(\mathfrak{N}\circ\mathfrak{F})}(G)/Z_{\pi\mathfrak{N}}(G^\mathfrak{F})=Z_{\pi\mathfrak{F}}(G/Z_{\pi\mathfrak{N}}(G^\mathfrak{F}))$.\par
\medskip
\noindent\textit{Proof.} Obviously, $Z_{\pi\mathfrak{N}}(G^\mathfrak{F}/Z_{\pi\mathfrak{N}}(G^\mathfrak{F}))=1$. By induction, Lemma 2.5(4) and Lemma 2.8(1), we may assume that $Z_{\pi\mathfrak{N}}(G^\mathfrak{F})=1$. Then by Lemma 2.12(2), $Z_{\pi(\mathfrak{N}\circ\mathfrak{F})}(G)\cap G^\mathfrak{F}=1$, and thus $Z_{\pi(\mathfrak{N}\circ\mathfrak{F})}(G)\leq C_G(G^\mathfrak{F})$. Since $\mathfrak{F}=\mathtt{S}\mathfrak{F}$, $H^\mathfrak{F}\leq G^\mathfrak{F}$ for every subgroup $H$ of $G$, and so $C_G(G^\mathfrak{F})\leq \mathcal{N}_{\pi'\mathfrak{F}}(G)$. Hence $Z_{\pi(\mathfrak{N}\circ\mathfrak{F})}(G)\leq \mathcal{N}_{\pi'\mathfrak{F}}(G)$. By Theorem B(1) and Theorem D, $\mathcal{N}_{\pi'\mathfrak{F}}^{\infty}(G)=Z_{\pi(\mathfrak{N}\circ\mathfrak{F})}(G)$. This implies that $\mathcal{N}_{\pi'\mathfrak{F}}^{\infty}(G)=\mathcal{N}_{\pi'\mathfrak{F}}(G)=Z_{\pi(\mathfrak{N}\circ\mathfrak{F})}(G)$. Suppose further that $\mathfrak{F}$ is saturated. Then by Lemma 2.12(3), $Z_{\pi(\mathfrak{N}\circ\mathfrak{F})}(G)=Z_{\pi\mathfrak{F}}(G)$, and so $\mathcal{N}_{\pi'\mathfrak{F}}^{\infty}(G)=\mathcal{N}_{\pi'\mathfrak{F}}(G)=Z_{\pi(\mathfrak{N}\circ\mathfrak{F})}(G)=Z_{\pi\mathfrak{F}}(G)$.\par
\medskip
\noindent\textbf{Lemma 3.5.} Let $\mathfrak{F}$ be a formation such that $\mathfrak{F}=\mathtt{S}\mathfrak{F}$. Then $\mathcal{N}_{\pi'\mathfrak{F}}^{\infty}(G)\leq \mbox{Int}_{\mathfrak{N}_\pi\circ\mathfrak{F}}(G)$ if one of the following holds:\par
(i) $G\in \mathfrak{S}_\pi\circ\mathfrak{F}$.\par
(ii) $\mathfrak{F}$ satisfies the $\pi$-boundary condition (I).\par
\medskip
\noindent\textit{Proof.} Let $H$ be any subgroup of $G$ such that $H\in \mathfrak{N}_\pi\circ\mathfrak{F}$. Then we only need to prove that $H\mathcal{N}_{\pi'\mathfrak{F}}^{\infty}(G)\in \mathfrak{N}_\pi\circ\mathfrak{F}$. By Lemma 2.5(2), we have that $\mathcal{N}_{\pi'\mathfrak{F}}^{\infty}(G)\leq \mathcal{N}_{\pi'\mathfrak{F}}^{\infty}(H\mathcal{N}_{\pi'\mathfrak{F}}^{\infty}(G))$. It follows that $H\mathcal{N}_{\pi'\mathfrak{F}}^{\infty}(G)/\mathcal{N}_{\pi'\mathfrak{F}}^{\infty}(H\mathcal{N}_{\pi'\mathfrak{F}}^{\infty}(G))\cong(H\mathcal{N}_{\pi'\mathfrak{F}}^{\infty}(G)/\mathcal{N}_{\pi'\mathfrak{F}}^{\infty}(G))/(\mathcal{N}_{\pi'\mathfrak{F}}^{\infty}(H\mathcal{N}_{\pi'\mathfrak{F}}^{\infty}(G))/\mathcal{N}_{\pi'\mathfrak{F}}^{\infty}(G))\in \mathfrak{N}_\pi\circ\mathfrak{F}$. Then by Corollary 3.2, $H\mathcal{N}_{\pi'\mathfrak{F}}^{\infty}(G)\in \mathfrak{N}_\pi\circ\mathfrak{F}$. Hence $\mathcal{N}_{\pi'\mathfrak{F}}^{\infty}(G)\leq \mbox{Int}_{\mathfrak{N}_\pi\circ\mathfrak{F}}(G)$.\par
\medskip
\noindent\textbf{Proof of Theorem C.} By Lemma 2.11(2), the canonical local definition $F$ of $\mathfrak{N}_\pi\circ\mathfrak{F}$ can be defined as follows: $F(p)=\mathfrak{G}_p\circ\mathfrak{F}$ for all $p\in \pi$; $F(p)=\mathfrak{N}_\pi\circ\mathfrak{F}$ for all $p\in \pi'$. Note that $Z_{\pi(\mathfrak{N}\circ\mathfrak{F})}(G)=Z_{\pi(\mathfrak{N}_\pi\circ\mathfrak{F})}(G)$ by Lemma 2.8(6). Then by \cite[Theorem A]{Guo}, (2) is equivalent to (3).\par
Next we show that (1) is equivalent to (3). Suppose that (3) holds, that is, $\mathfrak{F}$ satisfies the $\pi$-boundary condition (III). Then clearly, $\mathfrak{F}$ satisfies the $\pi$-boundary condition (I). Therefore, for every group $G$, we have that $\mathcal{N}_{\pi'\mathfrak{F}}^{\infty}(G)\leq \mbox{Int}_{\mathfrak{N}_\pi\circ\mathfrak{F}}(G)$ by Lemma 3.5. Since (2) is equivalent to (3), $\mbox{Int}_{\mathfrak{N}_\pi\circ\mathfrak{F}}(G)=Z_{\pi(\mathfrak{N}\circ\mathfrak{F})}(G)\leq \mathcal{N}_{\pi'\mathfrak{F}}^{\infty}(G)$ by Lemma 3.3. Consequently, $\mathcal{N}_{\pi'\mathfrak{F}}^{\infty}(G)=\mbox{Int}_{\mathfrak{N}_\pi\circ\mathfrak{F}}(G)$ holds for every group $G$, and so (3) implies (1).\par
Now suppose that $\mathcal{N}_{\pi'\mathfrak{F}}^{\infty}(G)=\mbox{Int}_{\mathfrak{N}_\pi\circ\mathfrak{F}}(G)$ holds for every group $G$, and there exists a prime $p\in \pi$ such that $\mbox{Crit}_\mathtt{S}(\mathfrak{G}_p\circ\mathfrak{F})\nsubseteq \mathfrak{N}_\pi\circ\mathfrak{F}$. Let $L$ be a group of minimal order in the set of all groups $G\in \mbox{Crit}_\mathtt{S}(\mathfrak{G}_p\circ\mathfrak{F})\backslash(\mathfrak{N}_\pi\circ\mathfrak{F})$. Then by Lemma 2.13, $L^{\mathfrak{N}_\pi\circ\mathfrak{F}}$ is the unique minimal normal subgroup of $L$ and $O_p(L)=\Phi(L)=1$. Hence by \cite[Chap. B, Theorem 10.3]{Doe}, there exists a simple $\mathbb{F}_pL$-module $P$ which is faithful for $L$. Let $V=P\rtimes L$. For any subgroup $H$ of $V$ such that $H\in \mathfrak{N}_\pi\circ\mathfrak{F}$, if $PH=V$, then $P\cap H\unlhd V$. This implies that $P\cap H=1$ for $P$ is a simple $\mathbb{F}_pL$-module, and so $H\cong V/P\cong L\notin \mathfrak{N}_\pi\circ\mathfrak{F}$, a contradiction. Hence $PH<V$. Then clearly, $PH\cap L<L$, and thus $PH/P=P(PH\cap L)/P\cong PH\cap L\in \mathfrak{G}_p\circ\mathfrak{F}$. It follows that $PH\in \mathfrak{G}_p\circ\mathfrak{F}\subseteq \mathfrak{N}_\pi\circ\mathfrak{F}$. Therefore, $P\leq \mbox{Int}_{\mathfrak{N}_\pi\circ\mathfrak{F}}(V)=\mathcal{N}_{\pi'\mathfrak{F}}^{\infty}(V)$. If $P\nleq \mathcal{N}_{\pi'\mathfrak{F}}(V)$, then $P\cap \mathcal{N}_{\pi'\mathfrak{F}}(V)=1$. Since $C_V(P)=P$, $\mathcal{N}_{\pi'\mathfrak{F}}(V)\leq C_V(P)=P$, and so $\mathcal{N}_{\pi'\mathfrak{F}}(V)=1$, which is absurd. Hence $P\leq \mathcal{N}_{\pi'\mathfrak{F}}(V)$. It follows that $P\leq N_V(L^\mathfrak{F}O_{\pi'}(V))$, and thereby $L^\mathfrak{F}O_{\pi'}(V)\unlhd V$. Then $L^\mathfrak{F}O_{\pi'}(V)\leq C_V(P)=P$. This induces that $L^\mathfrak{F}=1$. Therefore, $L\in \mathfrak{F}$, a contradiction. This shows that (1) implies (3). Consequently, (1) is equivalent to (3). The theorem is thus proved.\par
\medskip
\noindent\textbf{Proof of Theorem E.} We can prove the theorem similarly as in the proof of Theorem C by using \cite[Theorem 4.22]{Guo}.\par
\medskip
Now we give some conditions under which the formations satisfy the $\mathbb{P}$-boundary condition (I) (resp. the $\mathbb{P}$-boundary condition (II), the $\mathbb{P}$-boundary condition (III), the $\mathbb{P}$-boundary condition (III) in $\mathfrak{S}$). Recall that if $\sigma$ denotes a linear ordering on $\mathbb{P}$, then a group $G$ is called a Sylow tower group of complexion (or type) $\sigma$ if there exists a series of normal subgroups of $G$: $1=G_0\leq G_1\leq\cdots\leq G_n=G$ such that $G_i/G_{i-1}$ is a Sylow $p_i$-subgroup of $G/G_{i-1}$ for $1\leq i\leq n$, where $p_1\prec p_2\prec\cdots\prec p_n$ is the ordering induced by $\sigma$ on the distinct prime divisors of $|G|$. Let $\mathfrak{T}_\sigma$ denote the class of all Sylow tower groups of complexion $\sigma$. By \cite[Chap. IV, Example 3.4(g)]{Doe}, $\mathfrak{T}_\sigma$ is a saturated formation. Also, a formation $\mathfrak{F}$ is said to be a $\check{S}$-formation (or have the Shemetkov property) if $\mbox{Crit}_\mathtt{S}(\mathfrak{F})\subseteq \mbox{Crit}_\mathtt{S}(\mathfrak{N})\cup \{\mbox{cyclic groups of prime order}\}$. Clearly, $\mathfrak{N}_\pi$ is a $\check{S}$-formation. For details and more examples, see \cite[Section 3.5]{Guo4}. Moreover, a group $G$ is said to be $\pi$-closed if $G$ has a normal Hall $\pi$-subgroup. Let $\mathfrak{C}_{\pi}$ denote the formation of all $\pi$-closed groups.\par
\medskip
\noindent\textbf{Proposition 3.6.} A formation $\mathfrak{F}$ satisfies the $\mathbb{P}$-boundary condition (I) if one of the following holds:\par
(1) $\mathfrak{F}\subseteq \mathfrak{T}_\sigma$.\par
(2) $\mathfrak{F}$ is a $\check{S}$-formation.\par
(3) $\mathfrak{F}\subseteq \mathfrak{C}_{2}$.\par
(4) $\mathfrak{F}\subseteq \mathfrak{N}_{2}$.\par
\medskip
\noindent\textit{Proof.} (1) By \cite[Theorem 8]{Ros}, $\mbox{Crit}_\mathtt{S}(\mathfrak{T}_\sigma)\subseteq \mathfrak{S}$, and so $\mbox{Crit}_\mathtt{S}(\mathfrak{F})\subseteq \mathfrak{T}_\sigma\cup \mbox{Crit}_\mathtt{S}(\mathfrak{T}_\sigma)\subseteq \mathfrak{S}$.\par
Statement (2) is clear by definition.\par
(3) Note that $\mathfrak{C}_{2}$ is a $\check{S}$-formation by \cite[Remark]{Sta}. Then by Feit-Thompson Theorem, $\mbox{Crit}_\mathtt{S}(\mathfrak{F})\subseteq \mathfrak{C}_{2}\cup \mbox{Crit}_\mathtt{S}(\mathfrak{C}_{2})\subseteq \mathfrak{S}$.\par
The proof of statement (4) is similar to (3).\par
\medskip
\noindent\textbf{Proposition 3.7.} A formation $\mathfrak{F}$ satisfies the $\mathbb{P}$-boundary condition (II) if one of the following holds:\par
(1) $\mathfrak{F}\subseteq \mathfrak{N}$.\par
(2) $\mathfrak{F}\subseteq \mathfrak{G}_{2'}$ (equivalently, $2\notin\pi(\mathfrak{F})$).\par
\medskip
\noindent\textit{Proof.} (1) By \cite[Chap. IV, Satz 5.4]{Hup}, for any $p\in\mathbb{P}$, $\mbox{Crit}_\mathtt{S}(\mathfrak{G}_p\circ \mathfrak{N})=\mbox{Crit}_\mathtt{S}(\mathfrak{N}_{p'})\subseteq \mbox{Crit}_\mathtt{S}(\mathfrak{N})\subseteq\mathfrak{S}$. Hence $\mbox{Crit}_\mathtt{S}(\mathfrak{G}_p\circ \mathfrak{F})\subseteq \mathfrak{N}_{p'}\cup\mbox{Crit}_\mathtt{S}(\mathfrak{N}_{p'})\subseteq \mathfrak{S}$.\par
(2) Note that by \cite[Remark]{Sta}, $\mbox{Crit}_\mathtt{S}(\mathfrak{G}_2\circ \mathfrak{G}_{2'})=\mbox{Crit}_\mathtt{S}(\mathfrak{C}_{2})\subseteq \mathfrak{S}$, and for any odd prime $p$, $\mbox{Crit}_\mathtt{S}(\mathfrak{G}_p\circ \mathfrak{G}_{2'})=\mbox{Crit}_\mathtt{S}(\mathfrak{G}_{2'})=\{\mbox{cyclic group of order 2}\}\subseteq \mathfrak{S}$. Hence for any $p\in\mathbb{P}$, $\mbox{Crit}_\mathtt{S}(\mathfrak{G}_p\circ \mathfrak{F})\subseteq (\mathfrak{G}_p\circ \mathfrak{G}_{2'})\cup\mbox{Crit}_\mathtt{S}(\mathfrak{G}_p\circ \mathfrak{G}_{2'})\subseteq \mathfrak{S}$.\par
\medskip
Recall that a group $G$ is called $p$-decomposable if there exists a subgroup $H$ of $G$ such that $G=P\times H$ for some Sylow $p$-subgroup $P$ of $G$. Also, we use $\mathfrak{N}^r$ to denote the class of all groups $G$ with $l(G)\leq r$, where $l(G)$ is the Fitting length of $G$.\par
\medskip
\noindent\textbf{Proposition 3.8.} (1) Let $\mathfrak{F}$ be a formation with $\pi(\mathfrak{F})=\mathbb{P}$ such that $\mathfrak{F}\subseteq \mathfrak{N}$. Then $\mathfrak{F}$ satisfies the $\mathbb{P}$-boundary condition (III).\par
(2) Let $\mathfrak{L}$ be the formation of all $p$-decomposable groups. Then $\mathfrak{N}^r\circ\mathfrak{L}$ satisfies the $\mathbb{P}$-boundary condition (III) in $\mathfrak{S}$.\par
(3) Let $\mathfrak{F}$ be a formation with $\pi(\mathfrak{F})=\mathbb{P}$ such that $\mathfrak{F}\subseteq \mathfrak{N}$. Then $\mathfrak{N}^r\circ\mathfrak{F}$ satisfies the $\mathbb{P}$-boundary condition (III) in $\mathfrak{S}$.\par
\medskip
\noindent\textit{Proof.} Statement (1) was proved in \cite[Proposition 4.9(ii)]{Guo}, and statement (2) follows from \cite[Lemma 5.2]{Ski1} and \cite[Proposition 4.9(i)]{Guo}.\par
(3) By \cite[Chap. IV, Theorem 1.16]{Doe}, we have that $\mathfrak{F}=\mathtt{S}\mathfrak{F}$. It follows from (1) and \cite[Proposition 4.9(i)]{Guo} that $\mathfrak{N}^r\circ\mathfrak{F}$ satisfies the $\mathbb{P}$-boundary condition (III) in $\mathfrak{S}$.\par
\section{Applications}
In this section, we investigate the structure of groups $G$ whose minimal subgroups are contained in $\mathcal{N}_{\pi'\mathfrak{F}}^{\infty}(G)$. Let $\Psi_{p}(G)=\langle x|x\in G, o(x)=p\rangle$ if $p$ is odd, and $\Psi_2(G)=\langle x|x\in G, o(x)=2\rangle$ if the Sylow 2-groups of $G$ are quaternion-free, otherwise $\Psi_2(G)=\langle x|x\in G, o(x)=2\mbox{ or }4\rangle$.\par
\medskip
\noindent\textbf{Lemma 4.1.} Suppose that $\Psi_{p}(G^{\mathfrak{N}_p})\leq Z_{\mathfrak{N}_p}(G)$. Then $G\in \mathfrak{N}_p$.\par
\medskip
\noindent\textit{Proof.} By Lemma 2.8(3), for any subgroup $H$ of $G$, $\Psi_{p}(H^{\mathfrak{N}_p})\leq H\cap Z_{\mathfrak{N}_p}(G)\leq Z_{\mathfrak{N}_p}(H)$. Then by induction, $H\in \mathfrak{N}_p$. We may, therefore, assume that $G\in \mbox{Crit}_\mathtt{S}(\mathfrak{N}_p)$. It follows from \cite[Chap. IV, Satz 5.4]{Hup} that $G^{\mathfrak{N}_p}$ is a Sylow $p$-subgroup of $G$. By \cite[Theorem 1.1]{Sem}, $G^{\mathfrak{N}_p}/\Phi(G^{\mathfrak{N}_p})$ is a $G$-chief factor, and the exponent of $G^{\mathfrak{N}_p}$ is $p$ or $4$ (when $p=2$ and $G^{\mathfrak{N}_p}$ is non-abelian). If $p=2$ and $G^{\mathfrak{N}_2}$ is non-abelian and quaternion-free, then by \cite[Theorem 3.1]{War}, $G^{\mathfrak{N}_2}$ has a characteristic subgroup $L$ of index 2. This induces that $L=\Phi(G^{\mathfrak{N}_2})$, and so $G^{\mathfrak{N}_2}$ is cyclic, which is contrary to our assumption. Hence $p$ is odd or $p=2$ and $G^{\mathfrak{N}_2}$ is either abelian or not quaternion-free. This implies that $G^{\mathfrak{N}_p}=\Psi_{p}(G^{\mathfrak{N}_p})\leq Z_{\mathfrak{N}_p}(G)$, and thereby $G\in \mathfrak{N}_p$. The lemma is thus proved.\par
\medskip
\noindent\textbf{Lemma 4.2.} Let $\mathfrak{F}$ be a saturated formation such that $\mathfrak{F}=\mathtt{S}\mathfrak{F}$ and $\pi\subseteq \pi(\mathfrak{F})$. Suppose that $\Psi_{p}(G^{\mathfrak{F}})\leq Z_{\pi\mathfrak{F}}(G)$ for every $p\in \pi$. Then $G\in \mathfrak{G}_{\pi'}\circ\mathfrak{F}$.\par
\medskip
\noindent\textit{Proof.} Assume that the result is false and let $G$ be a counterexample of minimal order. If $O_{\pi'}(G)>1$, then for every $p\in \pi$, $\Psi_{p}(G^{\mathfrak{F}}O_{\pi'}(G)/O_{\pi'}(G))\leq Z_{\pi\mathfrak{F}}(G)/O_{\pi'}(G)=Z_{\pi\mathfrak{F}}(G/O_{\pi'}(G))$ by Lemma 2.8(1). Hence by the choice of $G$, $G/O_{\pi'}(G)\in \mathfrak{G}_{\pi'}\circ\mathfrak{F}$, and thereby $G\in \mathfrak{G}_{\pi'}\circ\mathfrak{F}$, which is impossible. Therefore, $O_{\pi'}(G)=1$. Let $M$ be any maximal subgroup of $G$. Since $M^{\mathfrak{F}}\leq G^{\mathfrak{F}}$ by Lemma 2.1(2), for every $p\in \pi$, $\Psi_{p}(M^{\mathfrak{F}})\leq M\cap Z_{\pi\mathfrak{F}}(G)\leq Z_{\pi\mathfrak{F}}(M)$ by Lemma 2.8(3). Then by the choice of $G$, $M\in \mathfrak{G}_{\pi'}\circ\mathfrak{F}$. We may, therefore, assume that $G\in \mbox{Crit}_\mathtt{S}(\mathfrak{G}_{\pi'}\circ\mathfrak{F})$.\par
If $Z_{\pi\mathfrak{F}}(G)\nleq \Phi(G)$, then $G$ has a maximal subgroup $M$ such that $G=Z_{\pi\mathfrak{F}}(G)M$. It follows that $G/Z_{\pi\mathfrak{F}}(G)\in \mathfrak{G}_{\pi'}\circ\mathfrak{F}$. By Lemmas 2.8(4) and 2.8(6), $G\in \mathfrak{G}_{\pi'}\circ\mathfrak{F}$, a contradiction. Hence $Z_{\pi\mathfrak{F}}(G)\leq \Phi(G)$ is nilpotent. Since $O_{\pi'}(G)=1$, $Z_{\pi\mathfrak{F}}(G)$ is a $\pi$-group. Then it is easy to see that $Z_{\pi\mathfrak{F}}(G)=Z_{\mathfrak{F}}(G)$. By \cite[Chap. IV, Theorem 6.10]{Doe}, for every $p\in \pi$, $\Psi_{p}(G^{\mathfrak{F}})\leq Z_{\mathfrak{F}}(G)\cap G^{\mathfrak{F}}\leq Z(G^{\mathfrak{F}})$. It follows from Lemma 4.1 that $G^{\mathfrak{F}}\in \mathfrak{N}_\pi$. As $O_{\pi'}(G)=1$, we have that $G^{\mathfrak{F}}\in \mathfrak{N}\cap \mathfrak{G}_\pi$. Since $G\in \mbox{Crit}_\mathtt{S}(\mathfrak{G}_{\pi'}\circ\mathfrak{F})$ and $\mathfrak{F}=\mathtt{S}\mathfrak{F}$, $G\in \mbox{Crit}_\mathtt{S}(\mathfrak{F})$. Then a similar discussion as in the proof of Lemma 4.1 shows that $G^{\mathfrak{F}}$ is a $p$-group with $p\in\pi$ such that the exponent of $G^{\mathfrak{F}}$ is $p$ or $4$ (when $p=2$ and $G^{\mathfrak{F}}$ is not quaternion-free) by using \cite[Theorem 1.1]{Sem}. This implies that $G^{\mathfrak{F}}=\Psi_{p}(G^{\mathfrak{F}})\leq Z_{\mathfrak{F}}(G)$, and so $G\in \mathfrak{F}$. The final contradiction ends the proof.\par
\medskip
\noindent\textbf{Theorem 4.3.} Let $\mathfrak{F}$ be a formation such that $\mathfrak{F}=\mathtt{S}\mathfrak{F}$. Suppose that $\Psi_{p}(G^{\mathfrak{N}_\pi\circ\mathfrak{F}})\leq \mathcal{N}_{\pi'\mathfrak{F}}^{\infty}(G)$ for every $p\in \pi$ and one of the following holds:\par
(i) $G\in \mathfrak{S_{\pi}}\circ\mathfrak{F}$.\par
(ii) $\mathfrak{F}$ satisfies the $\pi$-boundary condition (II).\par
(iii) $2\in\pi$ and $\mathfrak{F}$ satisfies the $\{2\}$-boundary condition (II).\par
(iv) $\{2,q\}'\subseteq\pi$, where $q$ is an odd prime, and $\mathfrak{F}$ satisfies the $\{2,q\}'$-boundary condition (II).\par
\noindent Then $G\in \mathfrak{N}_\pi\circ\mathfrak{F}$.\par
\medskip
\noindent\textit{Proof.} If either the condition (i) or the condition (ii) holds, then by Theorem B(1), Theorem D and Lemma 2.8(6), $\Psi_{p}(G^{\mathfrak{N}_\pi\circ\mathfrak{F}})\leq \mathcal{N}_{\pi'\mathfrak{F}}^{\infty}(G)=Z_{\pi(\mathfrak{N}\circ\mathfrak{F})}(G)=Z_{\pi(\mathfrak{N}_{\pi}\circ\mathfrak{F})}(G)$ for every $p\in \pi$. Hence in both cases, the theorem follows from Lemma 4.2.\par
Now suppose that the condition (iii) holds. Then it is easy to see that $\mathcal{N}_{\pi'\mathfrak{F}}(G)\leq \mathcal{N}_{2'\mathfrak{F}}(G)$ by definition, and so $\mathcal{N}_{\pi'\mathfrak{F}}^{\infty}(G)\leq \mathcal{N}_{2'\mathfrak{F}}^{\infty}(G)$. It follows that $\Psi_{2}(G^{\mathfrak{N}_2\circ\mathfrak{F}})\leq \Psi_{2}(G^{\mathfrak{N}_\pi\circ\mathfrak{F}})\leq\mathcal{N}_{\pi'\mathfrak{F}}^{\infty}(G)\leq \mathcal{N}_{2'\mathfrak{F}}^{\infty}(G)$. By applying the condition (ii), $G\in \mathfrak{N}_2\circ\mathfrak{F}\subseteq \mathfrak{S}\circ\mathfrak{F}$, and thereby the condition (i) holds. Therefore, $G\in \mathfrak{N}_\pi\circ\mathfrak{F}$.\par
Finally, we assume that the condition (iv) holds. Then for every $p\in \{2,q\}'$, $\Psi_{p}(G^{\mathfrak{N}_{\{2,q\}'}\circ\mathfrak{F}})\leq \Psi_{p}(G^{\mathfrak{N}_\pi\circ\mathfrak{F}})\leq\mathcal{N}_{\pi'\mathfrak{F}}^{\infty}(G)\leq \mathcal{N}_{\{2,q\}\mathfrak{F}}^{\infty}(G)$. By applying the condition (ii), $G\in \mathfrak{N}_{\{2,q\}'}\circ\mathfrak{F}\subseteq \mathfrak{S}\circ\mathfrak{F}$ by Burnside's $p^aq^b$-theorem, and so the condition (i) holds. Hence $G\in \mathfrak{N}_\pi\circ\mathfrak{F}$.\par
\medskip
The next two corollaries can be regarded as generalizations of \cite[Theorem 5.2]{She} and \cite[Theorem 5.3]{She}, respectively.\par
\medskip
\noindent\textbf{Corollary 4.4.} Let $\mathfrak{F}$ be a formation such that $\mathfrak{F}=\mathtt{S}\mathfrak{F}$ and $\mathfrak{F}\subseteq \mathfrak{U}$. Suppose that all cyclic subgroups of $G$ of odd prime order are contained in $\mathcal{N}_{\mathfrak{F}}^{\infty}(G)$. Then:\par
(1) $G\in \mathfrak{S}$.\par
(2) The $p$-length of $G$ is at most 2 for every odd prime $p$, and if $\mathfrak{F}\subseteq \mathfrak{N}$, then the $p$-length of $G$ is at most 1 for every odd prime $p$.\par
(3) The Fitting length of $G$ is bounded by 4, and if $\mathfrak{F}\subseteq \mathfrak{N}$, then the Fitting length of $G$ is bounded by 3.\par
\medskip
\noindent\textit{Proof.} By \cite[Chap. IV, Satz 5.4]{Hup}, $\mathfrak{N}_2$ satisfies the $2'$-boundary condition (II), and so $\mathfrak{F}$ also satisfies the $2'$-boundary condition (II) for $\mathfrak{F}\subseteq \mathfrak{U}\subseteq \mathfrak{N}_2$. Since $\Psi_{p}(G^{\mathfrak{N}_{2'}\circ\mathfrak{F}})\leq \mathcal{N}_{\mathfrak{F}}^{\infty}(G)\leq \mathcal{N}_{\{2\}\mathfrak{F}}^{\infty}(G)$ for every odd prime $p$, by Theorem 4.3, $G\in \mathfrak{N}_{2'}\circ\mathfrak{F}\subseteq \mathfrak{S}$. Hence for every odd prime $p$, $G^{\mathfrak{U}}\leq G^\mathfrak{F}\in \mathfrak{N}_{2'}\subseteq \mathfrak{N}_p$, and so the $p$-length of $G$ is at most 2 for every odd prime $p$. It is clear that $G^{\mathfrak{N}^2}\leq G^\mathfrak{U}\in \mathfrak{N}_{2'}\subseteq \mathfrak{N}^2$. This implies that the Fitting length of $G$ is bounded by 4. Now consider that $\mathfrak{F}\subseteq \mathfrak{N}$. The discussion is similar as above.\par
\medskip
\noindent\textbf{Corollary 4.5.} Let $\mathfrak{F}$ be a formation such that $\mathfrak{F}=\mathtt{S}\mathfrak{F}$ and $\mathfrak{F}\subseteq \mathfrak{U}$. Suppose that all cyclic subgroups of $G$ of order prime or 4 are contained in $\mathcal{N}_{\mathfrak{F}}^{\infty}(G)$. Then:\par
(1) $G\in \mathfrak{S}$.\par
(2) The $p$-length of $G$ is at most 2 for every $p\in\mathbb{P}$, and if $\mathfrak{F}\subseteq \mathfrak{N}$, then the $p$-length of $G$ is at most 1 for every $p\in\mathbb{P}$.\par
(3) The Fitting length of $G$ is bounded by 3, and if $\mathfrak{F}\subseteq \mathfrak{N}$, then the Fitting length of $G$ is bounded by 2.\par
\medskip
\noindent\textit{Proof.} The corollary can be proved similarly as in the proof of Corollary 4.4.\par
\medskip
\noindent\textbf{Acknowledgments.}\par
The authors are very grateful to the referee for his/her careful reading and helpful comments.\par
\bibliographystyle{plain}
\bibliography{expbib}
\end{document}